\documentclass[oneside,english]{amsart}
\usepackage[T1]{fontenc}
\usepackage[latin9]{inputenc}
\usepackage{amsthm}
\usepackage[pdftex]{graphicx}

\makeatletter
\numberwithin{equation}{section}
\numberwithin{figure}{section}

\makeatother

\usepackage{babel}
\begin{document}

\title[\tiny{Some Properties of the Derivatives on Sierpinski Gasket Type Fractals}]{Some Properties of the Derivatives on Sierpinski Gasket Type Fractals}
\author{Shiping Cao}
\address{School of Physics, Nanjing University, Nanjing, 210093, P.R. China.}
\curraddr{} \email{shipingcao@hotmail.com}
\thanks{}

\author{Hua Qiu$^*$}
\address{Department of Mathematics, Nanjing University, Nanjing, 210093, P. R. China.}
\curraddr{} \email{huaqiu@nju.edu.cn}
\thanks{$^*$ Corresponding author. }
\thanks{The research of the second author was supported by the Nature Science Foundation of China, Grant 11471157, and the Nature Science Foundation of Jiangsu Province of China, Grant
BK20131265.}

\subjclass[2000]{Primary 28A80.}

\keywords{derivatives on fractals, harmonic function, fractal Laplacian, tangent on fractals, Sierpinski gasket}

\date{}

\dedicatory{}
\begin{abstract}
In this paper, we focus on Strichartz's derivatives, a family of derivatives including the normal derivative, on p.c.f. (post critically finite) fractals, which are defined at vertex points in the graphs that approximate the fractal. We obtain a weak continuity property of the derivatives for  functions in the domain of the Laplacian. For a function with zero normal derivative at any fixed vertex, the derivatives, including the normal derivatives of the neighboring vertices will decay to zero. The optimal rates of approximations  are described and several non-trivial examples are provided to illustrate that our estimates are sharp. We also study the boundness property of derivatives for functions in the domain of the Laplacian. A necessary condition for a function having a weak tangent of order one at a vertex point is provided. Furthermore, we give a counter-example  of a conjecture of Strichartz on the existence of higher order weak tangents. 
\end{abstract}
\maketitle

\section{Introduction }

The theory of  analysis on fractals, analogous to that theory on manifolds, has been being well developed. The pioneering work is the analytic construction of the Laplacians, for a class of self-similar fractals that include the Sierpinski gasket as a typical example, developed by Kigami {[}Ki1-Ki6{]}, in which the Laplacians are defined as  renormalized limits of graph Laplacians. There are a lot of works in exploring some properties of these fractal Laplacians that are natural analogs of those for the usual Laplacian. For related works see {[}BK, BST, DSV, GKQS, IPRRS, KL, KSS, L, MT, QS, RS, T1, S1-S7{]} and references therein. Especially, there were several works in creating a calculus on fractals{[}BSSY, DRS, Ku, NSTY, S5, T2{]}. 

Since the fractal Laplacian acts as a differential operator with order greater than one, in analogy with the usual Laplacians on manifolds which are of second order, it is natural to make clear what is the first order derivative or gradient.  There are two approaches to find the gradient. One is to view the Dirichlet form as an integral of the inner product of gradients, see {[}Ku2, Ki3{]} for some works on this approach. However, it seems that this could not give direct information for a pointwise gradient. The other is to find the pointwise gradient directly. A paper of Teplyaev {[}T2{]} has made a satisfactory definition of the gradient at the general points in fractals and obtained some properties. For the vertex points case, Strichartz {[}S5{]} starting from the normal derivative, introducing some other derivatives at a vertex point $x$, together with which, he made up a (local) gradient $df(x)$. Please see {[}T2{]} to find a description of the relations in between these different definitions and results of Kigami, Kusuoka, Teplyaev and Strichartz.

In this paper, we continue to study the properties of Strichartz's derivatives at vertex points on fractals. 

We begin by assuming that the fractal $K$ is the invariant set of a finite iterated function system (i.f.s.) of contractive similarities in some Euclidean space $\mathbb{R}^n$. We denote these mappings $\{F_i\}_{i=1,\cdots, N}$. Then $K$ is the unique nonempty compact set satisfying
$$K=\bigcup F_i K.$$  We define $W_n$ as the collection of words $w=(w_1,\cdots,w_n)$ of the length $|w|=n$ from the alphabet $(1,2,\cdots,N)$ and write $F_w=F_{w_1}\circ F_{w_2}\circ\dots\circ F_{w_n}$. 
We call $F_w K$ a cell of level $n$.  

We use Strichartz's definition of the p.c.f. self-similar sets. $K$ is a \textit{post critically finite (p.c.f.) self-similar set} if $K$ is connected, and there is a finite set $V_0\subset K$ called the \textit{boundary}, such that the intersection of the sets $F_wK$ and $F_{w'}K$ is contained in the intersection of the boundary of these sets, $F_wV_0$ and $F_{w'}V_0$,  for any two different words $w$ and $w'$ with the same length. 

We denote $V_n=\bigcup_{w\in W_n}F_wV_0$ and $V_*=\bigcup_{n\geq 0} V_n$. A point $x\in V_*$ is called a \textit{junction vertex} if there are at least two different $w, w'\in W_n$ such that $x\in F_wK\cap F_{w'}K$. Otherwise we call $x$ a \textit{nonjunction vertex}.

We assume that a regular harmonic structure is given on a p.c.f. self-similar fractal $K$.  Thus there exists a self-similar \textit{Dirichlet form} $\mathcal{E}$ on $K$. It means for functions $f: K\rightarrow \mathbb{R}$, one has
$$
\mathcal{E}(f)=\sum_{j=1}^{N}r_j^{-1}\mathcal{E}(f\circ F_j)
$$
for some choice of \textit{ renormalization factors } $r_1,\cdots, r_N\in (0,1)$.
This quadratic form is obtained from the approximating of renormalized limit of $\mathcal{E}_m(f):=\mathcal{E}_m(f,f)$ on the $m$-level approximating graphs, where the $m$-level bilinear form is defined as
$$\mathcal{E}_m(f,g)=\sum_{|w|=m}r_w^{-1}\mathcal{E}_0(f\circ F_w,g\circ F_w),$$
with
$$\mathcal{E}_0(f,g)=\sum_{1\leq i<j\leq N_0}c_{ij}(f(v_i)-f(v_j))(g(v_i)-g(v_j)),$$
for some positive conductances $c_{i,j}$. Here $r_w=r_{w_1}\cdots r_{w_m}$. 

Let $\mathcal{H}_0$ denote the space of \textit{harmonic functions} on $K$ that minimize $\mathcal{E}_m$ at all levels for the given boundary values on $V_0$. Let $\mathcal{H}_m$ denote the space of continuous functions whose restrictions to each $F_wK$, for $|w|=m$, are harmonic (i.e., $h\circ F_w$ is harmonic).

The reader is referred to the books {[}Ki7{]} and {[}S8{]} for exact definitions, and any unexplained notations. 

Two additional assumptions are made, which are same as Strichartz did in {[}S5{]}.

\textbf{Hypothesis 1.1.}  \textit{(a) Each point $v_j, j=1,2,\cdots,N_0$ in the boundary set $V_0$ is the fixed point of a unique mapping in the i.f.s., which we denote $F_j$. Also, we assume that for any $F_i$ and $F_j$ in the i.f.s., $i\neq j$, the intersection $F_iK\cap F_jK$ consists of at most one point $x$ with $x=F_i v_m=F_j v_n$ for some points $v_m$ and $v_n$ in $V_0$.} 

\textit{(b) For each $v_j\in V_0$, let $M_j$ denote the $N_0\times N_0$ matrix that transforms the values $h|_{V_0}$ to $h|_{F_jV_0}$ for harmonic functions $h$, i.e.,
$$h(F_jv_k)=\sum_{l=1}^{N_0}(M_{j})_{kl}h(v_l).$$
We assume that each $M_j$ has a complete set of real left eigenvectors $\beta_{jk}$ with real nonzero eigenvalues $\lambda_{jk}$, i.e.,
$$\beta_{jk}M_j=\lambda_{jk}\beta_{jk},$$
where for each $j$ the eigenvalues $\lambda_{jk}$ are labeled in decreasing order of absolute value.
}

We will list some basic properties of the eigenvalues and eigenvectors of the matrixes $M_j$ in the next section. But here we only mention that the largest eigenvalue of $M_j$ is $\lambda_{j1}=1$, the second largest eigenvalue is $\lambda_{j2}=r_j$, the $j$-th renormalization factor of the harmonic structure, the eigenspace of the second eigenvalue $\lambda_{j2}$ is of one dimension, and we have $|\lambda_{jk}|<\lambda_{j2}$ for $k\geq 3$.

The following is the definition of Strichartz's derivatives at the boundary points. 

\textbf{Definition 1.2.} \textit{Let $f$ be a continuous function defined in a neighborhood of $v_j$. The derivatives $d_{jk}f(v_j)$ for $2\leq k\leq N_0$ are defined by the following limits, if they exist,
$$d_{jk}f(v_j)=\lim_{m\rightarrow\infty}\lambda_{jk}^{-m}\beta_{jk}f|_{F_j^mV_0}$$ where
$\beta_{jk}f|_{F_j^mV_0}$ is
$$\sum_{l=1}^{N_0}(\beta_{jk})_lf(F_j^mv_l).$$}

The derivative associated with $\beta_{j2}$ is just a multiple of the normal derivative at $v_j$. We could view other derivatives are of somewhat higher "order".
If $h$ is harmonic in a neighborhood of $v_j$, then all derivatives $d_{jk}$ exist and may be evaluated without taking the limit. See Lemma 3.3 in {[}S5{]}.

The above definition could be extended to all vertex points in $V_*$.  For a nonjunction vertex $x$, we suppose $n$ is the first value such that $x\in V_n$. Then there is a unique word $w$ of length $n$ such that $x=F_wv_j$ for some $1\leq j\leq N_0$.  We write $U_{m}(x)=F_wF_j^mK$, and call $\{U_m(x)\}_{m\geq 0}$ a standard system of neighborhoods of $x$. For a junction vertex $x$, by the Hypotheses 1.1(a), it is just an image under a mapping $F_w$ of a junction vertex in $V_1$. Let $J(x)$ denote the set of indices $j$ such that there exists $j'$ with $F_w^{-1}x=F_jv_{j'}$. Obviously, $\sharp J(x)\geq 2$. Suppose $n$ is the first value such that $x\in V_n$, then there exists a word $w$ of length $n-1$, such that $x=F_wF_jv_{j'}$ for all $j\in J(x)$. Write $U_m(x)=\bigcup_{j\in J(x)}F_wF_jF_{j'}^mK$, and call $\{U_{m}(x)\}_{m\geq 0}$ a standard system of neighborhoods of $x$.

 \textbf{Definition 1.3.} \textit{Let $f$ be a continuous function defined in a neighborhood of a vertex $x\in V_n\setminus V_{n-1}$.}
 
 \textit{
   (a) If $x=F_wv_j$ is a nonjunction vertex, then the derivatives $d_{jk}f(x)$ for $2\leq k\leq N_0 $ are defined by the following limits, if they exist,
\begin{equation}
   d_{jk}f(x)=\lim_{m\to \infty}r_w^{-1}\lambda_{jk}^{-m}\beta_{jk}f|_{F_wF_j^mV_0}.\label{nonj}
\end{equation}}
\textit{
   (b) If $x$ is a junction vertex, then the derivatives $d_{j'k}f(x)$ for $j\in J(x)$ and $2\leq k\leq N_0$ are defined by the following limits, if they exist,
$$   d_{j'k}f(x)=\lim_{m\to \infty}r_w^{-1}r_j^{-1}\lambda_{j'k}^{-m}\beta_{j'k}f|_{F_wF_jF_{j'}^mV_0}  .
$$
   Furthermore, the normal derivatives $d_{j'2}f(x)$ are said to satisfy the compatibility condition if
$$
   \sum_{j\in J(x)}d_{j'2}f(x)=0.
$$}

We write $df(x)$ for the collection of all derivatives defined here, and refer to it as the \textit{gradient} of $f$ at $x$. $f$ is called \textit{differentiable} at vertex $x$ if all the derivatives at $x$ exist and the compatibility condition holds if $x$ is a junction vertex. If $h$ is harmonic in a neighborhood of $x$, then $h$ is differentiable and all the derivatives may be evaluated without taking the limit. See Lemma 3.6 in {[}S5{]}.

 \textit{Remark.}  For higher "order" derivatives $d_{jk}$ or $d_{j'k} (3\leq k\leq N_0)$, one could check that we have two different scalings. Let $x=F_wv_j$ be a nonjunction vertex. Then for any word $u$, we have
 $$d_{jk}(f\circ F_u^{-1})(F_ux)=r_u^{-1}d_{jk}f(x),$$ and for any $m$, we have
 $$d_{jk}(f\circ F_wF_j^mF_w^{-1})(x)=\lambda_{jk}^m d_{jk}f(x).$$
 The junction vertex case is very similar. we omit it.
 
Let $\mu$ be a self-similar measure on $K$ with  weights $(\mu_1,\cdots,\mu_{N})$. It was proved in {[}S5{]} that for a function  $f\in dom (\Delta_\mu)$, the normal derivatives $d_{j2}f(x)$ and $d_{j'2}f(x)$ are uniformly bounded as $x$ varies over all  vertices. And for a harmonic function $h$ which take zero normal derivative at a vertex $x$, the normal derivatives of its neighboring vertices will decay to zero. See Theorem 4.3 in {[}S5{]}.

 In this paper, we will extend the boundness property to all derivatives, and the  weak continuity property to functions in the domain of the Laplacian for all derivatives. Moreover, we obtain the exact rate of approximations. We will prove the following three theorems.  These results answer the question post by Strichartz in {[}S5{]} positively. 

\textbf{Theorem 1.4.} \textit{Let $f\in dom(\Delta_{\mu})$. Then the normal derivatives of $f(x)$ are uniformly bounded as $x$ varies over all vertices. Furthermore, For any fixed nonjunction vertex $x=F_wv_j$ (or junction vertex $x=F_wF_jv_{j'}$), if $d_{j2}f(x)=0$ (or $d_{j'2}f(x)=0$), then  
 \begin{equation}
      d_{i2}f(y) (\text{or }  (d_{i'2}f(y))=
      \begin{cases}
      O(\mu_{j}^m),\qquad &\text{if $r_j\mu_{j}>|\lambda_{j3}|,$}\\
      O(m\mu_{j}^m),\quad &\text{if $r_j\mu_{j}=|\lambda_{j3}|,$}\\
      O((\lambda_{j3}r_j^{-1})^m),  &\text{if $r_j\mu_{j}<|\lambda_{j3}|,$}
      \end{cases}
      \end{equation}
      for all vertices $y\in U_m(x).$}
 
 \textbf{Theorem 1.5.} 
      \textit{(a) Let $h$ be a harmonic function. Then all the derivatives of $h(x)$ are uniformly bounded as $x$ varies over all vertices. }
      
      \textit{
       (b) Assume  $r_i\mu_i<|\lambda_{iN_0}|$ for $1\leq i\leq N_0$. Let $f\in dom(\Delta_\mu)$, then $f$ is differentiable at all vertices and all the derivatives of $f(x)$ are uniformly bounded.}     
      
\textbf{Theorem 1.6.} 
      \textit{(a) Let $h$ be a harmonic function, $x=F_wv_j$ be a nonjunction vertex (or $x=F_wF_jv_{j'}$ be a junction vertex) with zero normal derivative. Then for any vertices $y\in U_m(x)\setminus\{x\}$, we have
   \begin{equation}
   d_{ik}h(y) (\text{or } d_{i'k}h(y))=O((\lambda_{j3}r_j^{-1})^m).
   \end{equation}}
   \textit{
   (b) Assume  $r_i\mu_i<|\lambda_{iN_0}|$ for $1\leq i\leq N_0$. Let $f\in dom(\Delta_\mu)$, and $x$ be a vertex with zero normal derivative, then the above estimate still holds, with $f$ replaced by $h$. }
   
   Several non-trivial examples will be provided to illustrate that our estimates are optimal. There are some important fractals, including the \textit{Sierpinski gasket}, for which the condition $r_i\lambda_i<|\lambda_{iN_0}|$ does not hold. However, for these fractals, the results in Theorem 1.5 and 1.6 are still valid, provided we assume that the function $\Delta_\mu f$ satisfies an appropriate H\"{o}lder condition. 
   
   These results will be given in Section 3 and Section 4.
   
   We also study tangent in this paper.  As in {[}S5{]}, for a function $f$ differentiable at a vertex $x$, a \textit{weak tangent of order one} is defined as a harmonic function on $U_0(x)$, which assumes the same value and the values of derivatives at $x$ as those of $f$, denoted as $T_1^x(f)$ at $x$.

In Theorem 3.11 in {[}S5{]}, it is proved that for any function $f$ which is differentiable at a vertex $x$, let $h_m$ denote the harmonic function that assume the same values as $f$ at the boundary points of $U_m(x)$, extended to be harmonic on $U_0(x)$, then $h_m$ converges uniformly to $T_1^x(f)$ on $U_0(x)$ as $m\rightarrow\infty$. However, we will prove that it is not true provided that some reasonable assumptions on the harmonic structure or even the self-similar measure of the fractal be added. 
   
 If we assume that $\sharp V_0=3$ and all structures have full $D3$ symmetry, we could extend the definition of one order tangent to higher order. Here $D3$ symmetry means that all the structures are invariant under any homeomorphism of $K$. In this case, all $r'$s and $\mu'$s should be the same. Denote $\rho$ the value of $r_j\mu_j$ for $j=1,2,3$. Denote $\lambda_3$ the value of $\lambda_{j3}$ for  $j=1,2,3$ since they are also the same. Then for a vertex $x$ and a function $f$ defined in a neighborhood of $x$, an $n$-harmonic function $h$ is called a \textit{weak tangent of order $n$} if
 \begin{equation}
 (f-h)|_{\partial U_m(x)}=o((\rho^{n-1}r)^m),\label{wt1}
 \end{equation}
 and 
 \begin{equation}\label{wt2}
 (f-h-(f-h)\circ{g_x})|_{\partial U_m(x)}=o((\rho^{n-1}\lambda_3)^m), 
 \end{equation}   where
 \textit{$n$-harmonic functions} means those functions satisfy the equation $\Delta_\mu^nh=0$, and $g_x$ is a local point symmetry at $x$ with reasonable understanding(we omit the exact definition). 
 
 In {[}S5{]}, there is a conjecture, Conjecture 6.7, says that for a function $f\in dom(\Delta_\mu^{n-1})$, $f$ has a weak tangent of order $n$ at $x$ if and only if $d\Delta_{\mu}^{k}f(x)$ exists with compatibility conditions holding at $x$ for each $k\leq n-1$. It is true for $n=1$ since it is exact the definition of one order tangent. However, it is not true for $n\geq 2$. We will give a counter-example. 
 
 These results about tangents will be given in Section 5. 
 
  This paper can be regarded as a supplement of {[}S5{]}. Before ending of this section,  we mention a very useful result which could be obtained by an easy combination of the results in the appendix of {[}S5{]} and the results in the appendix of {[}T2{]}, says that, any function $f$ in the domain of the Laplacian satisfies an estimate 
  \begin{equation}
  |f(x)-f(y)|\leq cr_w\label{holder}
  \end{equation} for any $x,y\in F_wK$, where the constant $c$ is a multiple of $\Vert f\Vert_\infty+ \Vert\Delta_\mu f\Vert_\infty$.

\section{Basic results of the eigenvectors of $M_j$}

In this section, we will give some basic properties of the eigenvalues and eigenvectors of the transformation matrix $M_j$. Let $\{\lambda_{jk}\}_{1\leq k\leq N_0}$ be the set of eigenvalues labeled in decreasing order of absolute value. For each $\lambda_{jk}$, we denote $\beta_{jk}$ and $\alpha_{jk}$ the left and right eigenvectors of $\lambda_{jk}$ respectively. Additionally, we normalize that $\beta_{jk}\alpha_{jk}=1$. 

\textbf{Proposition 2.1.}  \textit{(a) The largest eigenvalue of $M_j$ is $\lambda_{j1}=1$. It has a right eigenvector $\alpha_{j1}=(1,\cdots,1)^t$, and a left eigenvector $\beta_{j1}$  with $(\beta_{j1})_l=\delta_{jl}$.}

\textit{(b) The second largest eigenvalue is $\lambda_{j2}=r_j<1$, the $j$-th renormalization factor of the harmonic structure. It has a left eigenvector $\beta_{j2}$ with $(\beta_{j2})_j=\sum_i c_{ij}$ and $(\beta_{j2})_l=-c_{lj}$ for $l\neq j$.}

\textit{(c) The eigenspace of $\lambda_{j2}$ is of one dimension and $|\lambda_{jk}|<\lambda_{j2}$ for $k\geq 3.$}

\textit{(d) $\beta_{jk}\alpha_{jl}=\delta_{kl}$ for $1\leq k,l\leq N_0$, where $\beta_{jk}\alpha_{jl}$ is $\sum_{s=1}^{N_0}(\beta_{jk})_s(\alpha_{jl})_s$.}

\textit{(e) For $k\geq 2$, $\sum_{l=1}^{N_0}(\beta_{jk})_l=0$ and $(\alpha_{jk})_j=0$.}

\textit{Proof.} One could find the proofs of (a), (b), (c) from {[}S5{]}. (d) is obvious. (e) follows from the combining of (a) and (d). $\Box$

Let $\{h_{jk}\}_{1\leq k\leq N_0}$ be a collection of harmonic functions on $K$, where each $h_{jk}$ assume values $\alpha_{jk}$ on $V_0$, i.e., $h_{jk}(v_l)=(\alpha_{jk})_l$ for each $l$. Obviously, $h_{j1}$ assumes constant value $1$ on $K$.

\textbf{Proposition 2.2.} \textit{(a) $h_{jk}|_{F_jV_0}=\lambda_{jk}h_{jk}|_{V_0},$ $d_{jk}h_{jl}(v_j)=\delta_{kl}$. }

\textit{(b) $h_{jk}(v_j)=0$ for $k\geq 2$.}

\textit{(c) $\{h_{jk}\}_{1\leq k\leq N_0}$ spans the space of harmonic functions on $K$. For any harmonic function $h$, it could be written into a linear combination that
\begin{equation}
h(\cdot)=h(v_j)+\sum_{k=2}^{N_0}d_{jk}h(v_j)h_{jk}(\cdot).
\end{equation}}

\textit{Proof.} (a) follows from the definition of $\alpha_{jk}$ and $\beta_{jk}$. (b) follows from Proposition 2.1(e). (c) is a corollary of (a) and (b). $\Box$

In the rest of this section, we will give some necessary and sufficient conditions for $(\beta_{jk})_j=0$ for all $k\geq 3$, which means that in this case the calculation of high "order" derivatives of a function $f$ at $v_j$ will not involve the value $f(v_j)$. This will be useful in Section 5.

\textbf{Proposition 2.3.}  \textit {The following three conditions are equivalent.}

\textit{(a) $(\beta_{jk})_j=0$ for all $k\geq 3$.}

\textit{(b) $(\alpha_{j2})_l=c(1-\delta_{jl})$ for all $l$,  where $c$ is a nonzero constant.}

\textit{(c) The $j$-th column of $M_j$ assumes the values that $(M_j)_{lj}=1-\lambda_{j2}+\lambda_{j2}\delta_{jl}.$}

\textit{Proof. } (a)$\Rightarrow$(b) Combining (a) and Proposition 2.1(e), we have that $\beta_{jk}, k\geq 3$ expand the linear space of dimension $N_0-2$ orthogonal to the constant vector and $\delta_{jl}$. Since $\beta_{jk}\alpha_{j2}=0, k\geq 3$, we conclude that
$$
   (\alpha_{j2})_l=s+t\delta_{jl}, l\geq1,
$$
for some constants $s$ and $t$. Moreover, by Proposition 2.1(e), $(\alpha_{j2})_j=0$. This determines that $s=-t$, which immediately yields (b).

(b)$\Rightarrow$(c)  Taking $\alpha_{j2}$ into the characteristic equation, we have 
$$M_j\alpha_{j2}=\lambda_{j2}\alpha_{j2},$$
which yields that 
$$\sum_{k\neq j}(M_j)_{lk}=\lambda_{j2}, \text{ for all } l\neq j.$$
Noticing that all row sums of $M_j$ are one and the $j$-th row of $M_j$ is $\delta_{kj}$, we then have
$$(M_j)_{lj}=1-\lambda_{j2} \text { for } l\neq j, \text { and } (M_j)_{jj}=1,$$ which is what (c) says.

(c)$\Rightarrow$(a)  For each $k\geq 3$, since $\beta_{jk}M_j=\lambda_{jk}\beta_{jk}$, by considering the $j$-th column of $M_j$, we have
$$\sum_{l\neq j}(1-\lambda_{j2})(\beta_{jk})_l+(\beta_{jk})_j=\lambda_{jk}(\beta_{jk})_j.$$
Combining the above formula with Proposition 2.1(e),  we obtain that $(\beta_{jk})_j=0$. Thus (a) holds. $\Box$

\textit{Remark.} In the $D3$ symmetry case, condition (c) automatically holds. Thus $(\beta_{j3})_j=0$, which means that the tangential derivative of a function $f$ at $v_j$ does not involve the value $f(v_j)$. 

Finally, we give an example which does not satisfy the conditions in Proposition 2.3.

\textit{Example 2.4.} 
Let $v_1, v_2, v_3$ be the vertices of an equilateral triangle and let $F_i(x)=\frac{1}{2}(x+v_i)$, i=1,2,3. The \textit{Sierpinski gasket}, $\mathcal{SG}$, is the unique compact set such that $\mathcal{SG}=\bigcup_{i=1}^3F_i\mathcal{SG}$. Then $V_0=\{v_1,v_2,v_3\}$. 

Consider a family of self-similar Dirichlet forms on $\mathcal{SG}$, that has a single bilateral symmetry. So we require $r_2=r_3$ and
$$\mathcal{E}_0(f)=(f(v_1)-f(v_2))^2+(f(v_1)-f(v_3))^2+c(f(v_2)-f(v_3))^2$$ 
for some $c>0$. We denote the conductances of $\mathcal{E}_0$ and $r_2\mathcal{E}_1$ on the edges  of the graphs $\Gamma_0$ and $\Gamma_1$ in Figure 2.1. where $s=r_2/r_1$ is a constant to be determined. 

\begin{figure}[h]
\begin{center}
\includegraphics[width=4.5cm,totalheight=4cm]{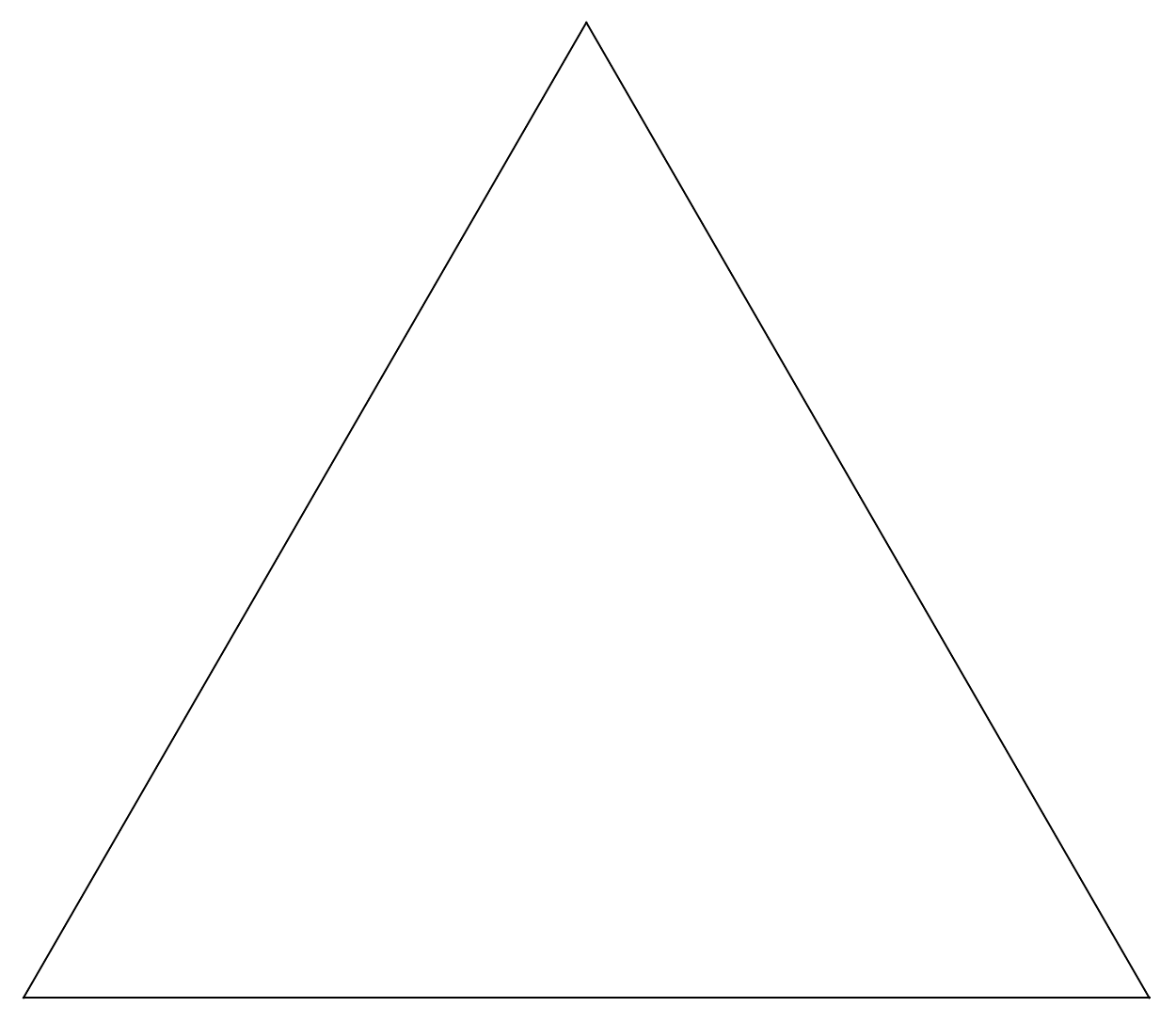}\hspace{2cm}
\includegraphics[width=4.5cm,totalheight=4cm]{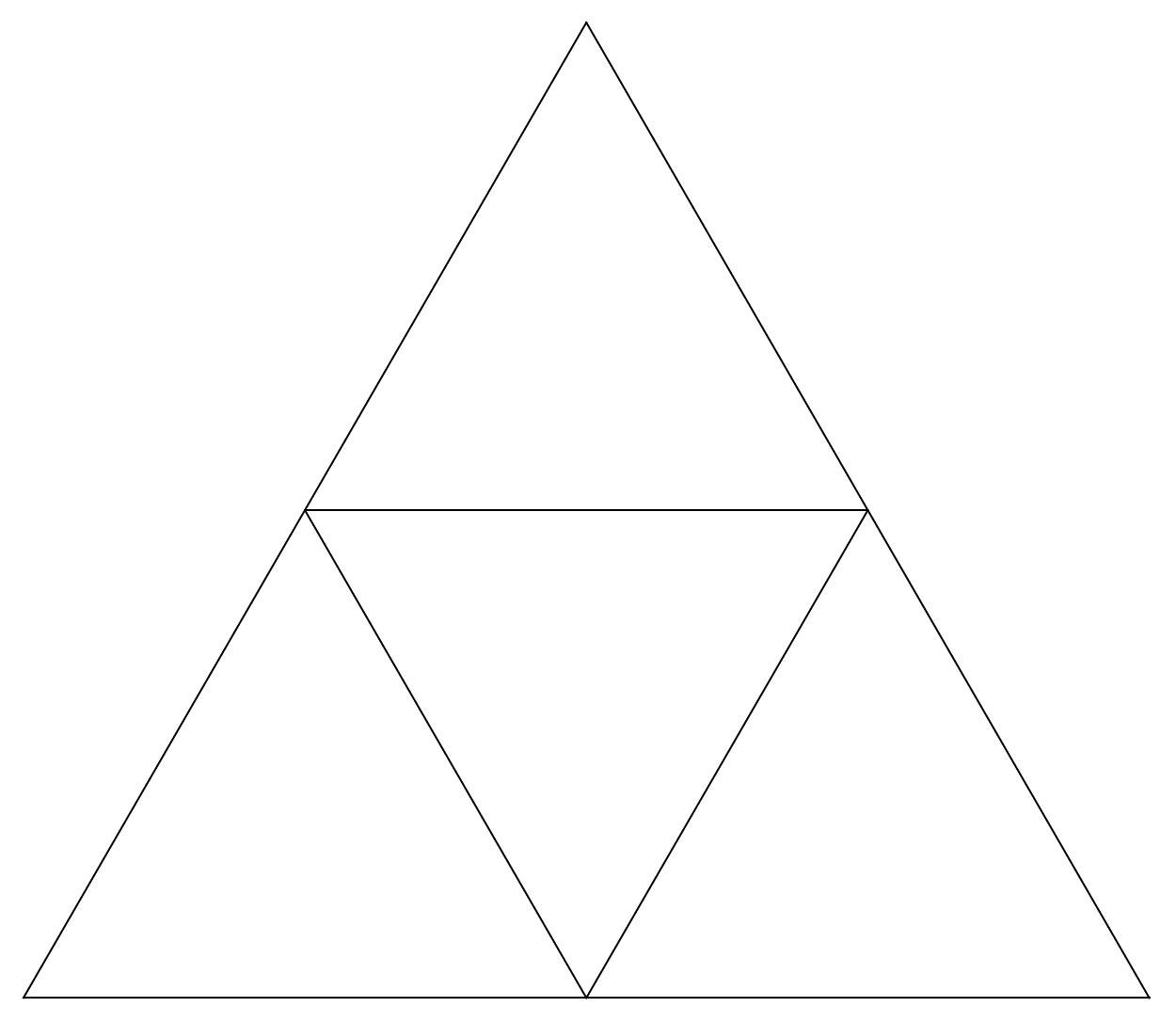}
\setlength{\unitlength}{1cm}
\begin{picture}(0,0) \thicklines
\put(-11.7,-0){$v_2$}
\put(-6.9,-0){$v_3$}
\put(-9.2,4.0){$v_1$}
\put(-10.5,2.0){$1$}
\put(-7.9,2.0){$1$}
\put(-9.2,0.2){$c$}
\put(-5.1,-0){$v_2$}
\put(-0.3,-0){$v_3$}
\put(-2.6,4.0){$v_1$}
\put(-3.3,2.9){$s$}
\put(-1.8,2.9){$s$}
\put(-2.7,2.1){$sc$}
\put(-4.4,0.9){$1$}
\put(-2.9,0.9){$1$}
\put(-3.65,0.2){$c$}
\put(-2.25,0.9){$1$}
\put(-0.75,0.9){$1$}
\put(-1.50,0.2){$c$}
\end{picture}
\begin{center}
\vspace{0.4cm}
\textbf{Figure 2.1. The conductances of $\mathcal{E}_0$ and $r_2\mathcal{E}_1$.}
\end{center}
\end{center}
\end{figure}

The renormalization equation requires $s$ and $c$ has the relationship
$$
         3s^2c^2+2s^2c-2sc^2-2c-1=0.
$$
         A detailed calculation could be found in Chapter 4 of {[}S8{]}.

Let $h$ be a harmonic function on $\mathcal{SG}$ with respect to the above Dirichlet form. The mean value equations of $h$ at vertices $F_2v_3, F_1v_3$ and $F_2v_1$ give that
$$
         \begin{cases}
         (2+2c)h(F_2v_3)-h(F_1v_3)-h(F_2v_1)-ch(v_2)-ch(v_3)=0,\\
         (2+s+sc)h(F_1v_3)-h(F_2v_3)-sch(F_2v_1)-sh(v_1)-h(v_3)=0,\\
         (2+s+sc)h(F_2v_1)-h(F_2v_3)-sch(F_1v_3)-sh(v_1)-h(v_2)=0.
         \end{cases}
 $$
This yields
         $$
         \begin{pmatrix} h(F_2v_3)\\h(F_1v_3)\\h(F_2v_1)\end{pmatrix}=
         \begin{pmatrix} 1-2\eta&\eta&\eta\\
                         \frac{1+s-2\eta}{2+s}&\frac{\eta}{2+s}+\frac{sc}{(2+s)(2+s+2sc)}&\frac{\eta}{2+s}+\frac{2+s+sc}{(2+s)(2+s+2sc)}\\
                         \frac{1+s-2\eta}{2+s}&\frac{\eta}{2+s}+\frac{2+s+sc}{(2+s)(2+s+2sc)}&\frac{\eta}{2+s}+\frac{sc}{(2+s)(2+s+2sc)}\end{pmatrix}
         \begin{pmatrix}h(v_1)\\h(v_2)\\h(v_3)\end{pmatrix}               
         $$
         where $\eta=\frac{2c+sc+1}{2sc+2s+4c+2}$. 
         
Thus the transformation matrix $M_2$ is 
$$
         M_2=\begin{pmatrix} \frac{1+s-2\eta}{2+s}&\frac{\eta}{2+s}+\frac{2+s+sc}{(2+s)(2+s+2sc)}&\frac{\eta}{2+s}+\frac{sc}{(2+s)(2+s+2sc)}\\
         0&1&0\\
         1-2\eta&\eta&\eta\end{pmatrix}.
$$
One can check that $M_2$ satisfies Hypothesis 1.1(b) when $|s-1|$ is sufficiently small. In fact, when $s=1$, $M_2$ is diagonalizable with three different eigenvalues and all entries of $M_2$ are continuous functions of $s$. 

Comparing $({M_2})_{12}$ and $({M_2})_{32}$, we can find they are not equal, since it leads to a different equation  
$$
         2s^2c^2+cs^2+cs-2c-s-1=0.
 $$
 
 Hence $M_2$ does not satisfies the condition(c) in Proposition 2.3, at least for those $s$ very close, but not equal to $1$, which means $(\beta_{23})_2\neq 0$.

\section{Boundness and weak continuity of normal derivatives}

We prove Theorem 1.4 in this section, and provide some examples to show that our results are sharp.

\textbf{Lemma 3.1.}  \textit{Let $f\in dom(\Delta_{\mu})$. Then the normal derivatives of $f$ are uniformly bounded on all vertices of $K$.}

This result is proved in {[}S5{]}, using Gauss-Green formula.  For the convenience of readers, we still provide a proof. But our proof is somewhat different to that in  {[}S5{]}, and could be extended to other derivatives. 

\textit{Proof.}  
Notice that from Proposition 2.1(e), for $1\leq j\leq N_0$, we have $\sum_{l=1}^{N_0}(\beta_{j2})_l=0$. Combine it with formula \eqref{holder}, the H\"{o}lder estimate of $f$, we obtain an estimate that
$$
      |r_w^{-1}\beta_{j2}f|_{F_wV_0}|\leq c
$$
      for all words $w$ and all $j$, with some constant $c>0$. Since we have the existences of normal derivatives at all vertices, we get that
   $$
      |d_{j2}f(x)|=|\lim_{m\to \infty}r_w^{-1}r_j^{-m}\beta_{j2}f|_{F_wF_j^mV_0}|\leq c. \quad\quad\Box
$$

Now we devote to prove the weak continuity. For convenience, we give the proof in the case of $x=v_j$,  since for other vertices, we only need to use scaling.  First, we give some lemmas.

\textbf{Lemma 3.2.} \textit{Let $h$ be any harmonic function with $d_{j2}h(v_j)=0$. Then for any vertices $y\in F_{j}^mK$, we have
$$
            d_{i2}h(y) (\text{or } d_{i'2}h(y))=O((\lambda_{j3}r_j^{-1})^m).
$$
}
\textit{Proof.} Since $h$ assumes $0$ normal derivative at $v_j$, by using Proposition 2.2(c), we have 
$$h(\cdot)=h(v_j)+\sum_{k=3}^{N_0}d_{jk}h(v_j)h_{jk}(\cdot).$$
 So we need to prove the lemma for each $h_{jk}, k\geq 3$. Let $c$ denote the upper bound of normal derivatives for all $h_{jk}$ as guaranteed by Lemma 3.1. Then for any $y\in F_j^mK$,
\begin{equation}\label{34}
      \begin{aligned}
      |d_{i2}h_{jk}(y)|&=|r_j^{-m}d_{i2}(h_{jk}\circ F_j^m)(F_j^{-m}y)|\\&=|(r_j^{-1}\lambda_{jk})^md_{i2}h_{jk}(F_j^{-m}y)|\leq c(r_j^{-1}|\lambda_{jk}|)^m\leq c(r_j^{-1}|\lambda_{j3}|)^m,
      \end{aligned}
   \end{equation}
      by using Proposition 2.2(a) and Lemma 3.1.
       $\Box$
       
We will need the local Green's function. Recall that If $G(x,z)$ denotes the Green's function for the Dirichlet problem on $K$, then $G(x,z)$ has the following expression. 
$$
      G(x,z)=\sum_{|w|\geq 0}r_w\Psi(F_w^{-1}x,F_w^{-1}z),
$$
      where the summation is over all words, and $\Psi$ is a linear combination of products $\psi_p(x)\psi_q(z)$ where $\psi_p$ are tent functions in $\mathcal{H}_1$, taking value $1$ at $p\in V_1\setminus V_0$ and $0$ at other vertices of $V_1$. For each term $\Psi(F_w^{-1}x, F_w^{-1}z)$, the understanding is that it assumes value $0$ unless $x$ and $z$ both belong to the cell $F_wK$.
See detailed explanations in {[}Ki8{]} and {[}S5{]}.

\textbf{Lemma 3.3.}   \textit{Let $g\in C(F_wK)$. Then
$$
      u(x)=\int_{F_wK}r_wG(F_w^{-1}x,F_w^{-1}z)g(z)d\mu(z)
$$ solves the local Dirichlet problem that $-\Delta_\mu u(x)=g(x)$ on $F_wK$ and $u|_{F_wV_0}=0$. Furthermore, for each boundary vertex $F_wv_i$, 
\begin{equation}
             \partial_n u(F_wv_i)=-\int_{F_wK}H_i(F_w^{-1}z)g(z)d\mu(z),\label{deri}
\end{equation}
where $H_i$ is the harmonic function on $K$ taking values $H_i(v_j)=\delta_{ij}$. 
}
    
 \textit{Proof. }    Rewrite the integral by scaling,
$$
      u(x)=r_w\mu_w\int_KG(F_w^{-1}x,z)g\circ F_w(z)d\mu(z).
$$
We then have
$$
      -\Delta_\mu (u\circ F_w)(F_w^{-1}x)=r_w\mu_w(g\circ F_w)(F_w^{-1}x).
$$
Combining it with the fact that $\Delta_\mu(u\circ F_w)=r_w\mu_w(\Delta_\mu u)\circ F_w$, we get that $-\Delta_\mu u(x)=g(x)$ on $F_wK$. The Dirichlet boundary condition can be checked directly. 
\eqref{deri} could be derived from the Gauss-Green formula since $u$ has the zero boundary condition on $F_wK$.  $\Box$

We also need to estimate the derivatives of the function $\Psi$.

For $1\leq i\leq N_0$, $2\leq k\leq N_0$, by the definition of the function $\Psi$, there exists a piecewise harmonic function $a_{ik}\in\mathcal{H}_1$ satisfying 
      \begin{equation}
      a_{ik}(z)=d_{ik}\Psi(\cdot,z)(v_i)
      \end{equation}
    Obviously, $a_{ik}|_{V_0}=0$.
    
 \textbf{Lemma 3.4.}  \textit{Let $u\in dom(\Delta_\mu)$. Then for all $1\leq i\leq N_0$, $2\leq k\leq N_0$, $m\geq 0$, 
 \begin{equation}
 \lambda_{ik}^{-m}\beta_{ik}u|_{F_i^mV_0}=\beta_{ik}u|_{V_0}-\sum_{n=0}^{m-1}r_i^n\lambda_{ik}^{-n}\int_K a_{ik}\circ F_{i}^{-n}(z)\Delta_\mu u(z)d\mu(z).\label{mderi}
 \end{equation}
 }  
 
 \textit{Proof.}  Let $h$ be a harmonic function which assumes the same values on $V_0$ as $u$. Then 
 $$u=-\int_KG(x,z)\Delta_{\mu}u(z)d\mu(z)+h.$$ Taking the above formula into the left side of \eqref{mderi}, we obtain that it equals to
 \begin{eqnarray*}
  &&\lambda_{ik}^{-m}\beta_{ik}h|_{F_i^mV_0}-\lambda_{ik}^{-m}\int_{K}\beta_{ik}G(\cdot,z)|_{F_i^mV_0}\Delta_{\mu}u(z)d\mu(z)\\
  &=&\beta_{ik}u|_{V_0}-\sum_{n=0}^{m-1}r_i^n\int_{K}\lambda_{ik}^{-m}\beta_{ik}\Psi(F_i^{-n}\cdot, F_i^{-n}z)|_{F_i^mV_0}\Delta_{\mu}u(z)d\mu(z)\\
  &=&\beta_{ik}u|_{V_0}-\sum_{n=0}^{m-1}r_i^{n}\lambda_{ik}^{-n}\int_{K}\lambda_{ik}^{-(m-n)}\beta_{ik}\Psi(\cdot, F_i^{-n}z)|_{F_i^{m-n}V_0}\Delta_{\mu}u(z)d\mu(z)\\
     &=&\beta_{ik}u|_{V_0}-\sum_{n=0}^{m-1}r_i^n\lambda_{ik}^{-n}\int_{K}a_{ik}(F_i^{-n}z)\Delta_{\mu} u(z)d\mu(z),
 \end{eqnarray*}
 where we use the fact that $h$ is harmonic, $h|_{V_0}=u|_{V_0}$ and $\Psi(x,z)$ is piecewise harmonic with respect to the first variable $x$. $\Box$
    
 \textbf{Lemma 3.5.}   \textit{$\sum_{n=0}^\infty a_{i2}\circ F_i^{-n}(x)=-H_i(x),$ for each point $x\in K\setminus\{v_i\}$ and $\sum_{n=0}^\infty a_{i2}\circ F_i^{-n}(v_i)=0$.}
 
  Here for each term, one should understand that $a_{i2}\circ{F_i^{-n}}(x)$ is zero unless $x\in F_i^nK$.
 Thus for $x\in F_i^{m-1}K\setminus F_i^{m}K,$ the above summation only involves the first $m$ nonzero terms.    
    
\textit{Proof.}   Let $u$ be a function in $dom(\Delta_\mu)$ that satisfies the  Dirichlet boundary condition. We have 
$$
      u(x)=-\int_K G(x,z)\Delta_\mu u(z)d\mu(z).
$$
 Using Lemma 3.4, noticing that $r_i=\lambda_{i2}$, we then have  
 $$\lambda_{i2}^{-m}\beta_{i2}u|_{F_i^mV_0}=-\sum_{n=0}^m \int_K a_{i2}\circ F_i^{-n}(z)\Delta_\mu u(z)d\mu(z).$$   
Observing that $\sum_{n=0}^\infty |a_{i2}\circ F_i^{-n}(z)\Delta_\mu u(z)|$ is integrable, letting  $m\rightarrow\infty$, using Lebesque's Control-Convergent theorem, we get
 $$ \partial_nu(v_i)=-\int_K \sum_{n=0}^\infty a_{i2}\circ F_i^{-n}(z)\Delta_\mu u(z)d\mu(z).$$
 
 On the other hand, by using the Gauss-Green formula, noticing that $u$ satisfies the Dirichlet boundary condition, we also get that $$\partial_nu(v_i)=\int_K H_i\Delta_\mu ud\mu.$$
 
 Thus we have 
$$
      \sum_{n=0}^\infty a_{i2}\circ F_i^{-n}(z)=-H_i(z),\quad \text{ a.e.-}\mu.
$$
from the arbitrariness  of function  $u$.

Moreover we have that the function $\sum_{n=0}^\infty a_{i2}\circ F_i^{-n}(z)$ is continuous on $K\setminus\{v_i\}$. Thus we can remove the requirement of "a.e." from the above formula for those $x\neq v_i$. As for $x=v_i$,  it is easy to check $\sum_{n=0}^\infty a_{i2}\circ F_i^{-n}(v_i)=0\neq H_i(v_i)$ since $a_{i2}$ satisfies the Dirichlet boundary condition. $\Box$

  \textit{Proof of Theorem 1.4.}  Let $f\in dom(\Delta_\mu)$. The boundness property of the normal derivatives of $f$ comes from Lemma 3.1. Thus we only need to prove the weak continuity property. 
  
  As stated before, we only give the proof in the case of $x=v_j$. Without loss of generality, we assume $y=F_j^mF_\tau v_i$. To study $d_{i2}f(y)$, we need to study the behavior of $f$ in the cell $F_j^mF_\tau K$.

   From Lemma 3.3, 
   $$\partial_n\left(-\int_KG(\cdot,z)\Delta_\mu f(z)d\mu(z)\right)(v_j)=\int_K H_j(z)\Delta_\mu f(z)d\mu(z).$$ 
   So 
    $$\partial_n\left(-\int_K(G(\cdot,z)+h_{j2}(\cdot)H_j(z))\Delta_\mu f(z)d\mu(z)\right)(v_j)=0,$$ 
    by using Proposition 2.2(a).
    
    Noticing that $\partial_n f(v_j)=0$, the difference between the function  $$-\int_K(G(\cdot,z)+h_{j2}(\cdot)H_j(z))\Delta_\mu f(z)d\mu(z)$$ and $f$ is a harmonic function with $0$ normal derivative at $v_j$, denoted by $h$. Thus we could write
\begin{equation}
      f(\cdot)=-\int_K (G(\cdot,z)+h_{j2}(\cdot)H_j(z))\Delta_\mu f(z)d\mu(z)+h(\cdot).\label{f}
\end{equation}

By Lemma 3.2, we only have to estimate the normal derivatives of the first summand  of the right side of \eqref{f}. 

Since we only interest in the values of $f$ at those points in $F_j^mK$, we can rewrite the Green's function in the above integral for variables in $F_j^mK$ as
$$
      \begin{aligned}
      G(\cdot,z)&=\sum_{n=0}^{m-1}r_j^n\Psi(F_j^{-n}\cdot,F_j^{-n}z) \\&+\sum_{|w|\geq0}r_j^mr_w\Psi(F_w^{-1}\circ F_j^{-m}\cdot,F_w^{-1}\circ F_j^{-m}z).
      \end{aligned}
$$
Taking this expression into \eqref{f}, we could write $$f=f_1+f_2+h\text{ on } F_j^mK,$$ where
$$
      f_1(\cdot)=-\int_K(\sum_{n=0}^{m-1}r_j^n\Psi(F_j^{-n}\cdot,F_j^{-n}z)+h_{j2}(\cdot)H_j(z))\Delta_\mu f(z)d\mu(z),$$ and
$$f_2(\cdot)=-\int_{F_j^mK}\sum_{|w|\geq 0}r_j^mr_w\Psi(F_w^{-1}\circ F_j^{-m}\cdot,F_w^{-1}\circ F_j^{-m}z)\Delta_\mu f(z)d\mu(z).
$$

We estimate the normal derivatives of these two functions separately. 

First, we deal with $f_1$. 
By using Lemma 3.5,  and Proposition 2.2(c), we could decompose the function $\sum_{n=0}^{m-1}r_j^n\Psi(F_j^{-n}\cdot,F_j^{-n}z)+h_{j2}(\cdot)H_j(z)$, which is harmonic  on  $F_j^mK$, as
$$
      \begin{aligned}
      &\sum_{n=0}^{m-1}r_j^n\Psi(F_j^{-n}\cdot,F_j^{-n}z)+h_{j2}(\cdot)H_j(z)\\=&\sum_{n=0}^{m-1}\sum_{k=3}^{N_0}r_j^n\lambda_{jk}^{-n}a_{jk}\circ F_j^{-n} (z)h_{jk}(\cdot)-\sum_{n=m}^{\infty}a_{j2}\circ F_j^{-n}(z)h_{j2}(\cdot),
      \end{aligned}
$$
for variables in $F_j^mK$.

So on $F_j^mK$, we have the exact formula for $f_1$, 
\begin{equation}
            \begin{aligned}
            f_1(\cdot)=&-\sum_{k=3}^{N_0}h_{jk}(\cdot)\int_K\sum_{n=0}^{m-1}r_j^n\lambda_{jk}^{-n}a_{jk}\circ F_j^{-n}(z)\Delta_\mu f(z)d\mu(z)\\&+
            h_{j2}(\cdot)\int_K \sum_{n=m}^{\infty} a_{j2}\circ F_j^{-n}(z)\Delta_\mu f(z)d\mu(z).\label{31}
            \end{aligned}
\end{equation}

For each coefficient $\int_K\sum_{n=0}^{m-1}r_j^n\lambda_{jk}^{-n}a_{jk}\circ F_j^{-n}(z)\Delta_\mu f(z)d\mu(z)$ of $h_{jk}$, $k\geq 3$, we have the estimate that
\begin{equation}
     |\int_K\sum_{n=0}^{m-1}r_j^n\lambda_{jk}^{-n}a_{jk}\circ F_j^{-n}(z)\Delta_\mu f(z)d\mu(z)|=
     \begin{cases}
     O(1),\qquad &\text{if } r_j\mu_{j}<|\lambda_{jk}|,\\
     O(m),\qquad &\text{if } r_j\mu_{j}=|\lambda_{jk}|,\\\label{32}
     O(\mu_j^mr_j^{m}\lambda_{jk}^{-m}), &\text{if } r_j\mu_{j}>|\lambda_{jk}|.
     \end{cases}
\end{equation}
For the coefficient $\int_K \sum_{n=m}^{\infty} a_{j2}\circ F_j^{-n}(z)\Delta_\mu f(z)d\mu(z)$ of $h_{j2}$, 
we have
\begin{equation}
|\int_K \sum_{n=m}^{\infty} a_{j2}\circ F_j^{-n}(z)\Delta_\mu f(z)d\mu(z)|=O(\mu_j^m).\label{33}
\end{equation}
Combining \eqref{31}, \eqref{32}, \eqref{33} and the estimates \eqref{34} in the proof of Lemma 3.2 for normal derivatives of $h_{jk}, k\geq 3$ over $F_j^mK$, we have 
   \begin{equation}\label{f1}
      d_{i2}f_1(y)=
      \begin{cases}
      O(\mu_{j}^m),\qquad &\text{if $r_j\mu_{j}>|\lambda_{j3}|,$}\\
      O(m\mu_{j}^m),\quad &\text{if $r_j\mu_{j}=|\lambda_{j3}|,$}\\
      O((\lambda_{j3}r_j^{-1})^m),  &\text{if $r_j\mu_{j}<|\lambda_{j3}|,$}
      \end{cases}
      \end{equation}   
     for any $y\in F_j^mK$.

 Next, we estimate the normal derivatives of $f_2$ on $F_j^mK$.  For $y=F_j^mF_\tau v_i$, we can further divide $f_2$ into two functions on $F_j^mK$ as
\begin{equation}\label{310}
f_3(\cdot)=-\sum_{0\leq n\leq |\tau|-1}\int_{K} r_j^mr_{\tau_1}...r_{\tau_n}\Psi(F_{\tau_n}^{-1}...F_{\tau_1}^{-1}\circ F_j^{-m}\cdot,F_{\tau_n}^{-1}...F_{\tau_1}^{-1}\circ F_j^{-m}z)\Delta_\mu f(z)d\mu(z),
\end{equation} and
\begin{equation}\label{311}
f_4(\cdot)=-\int_K\sum_{w}r_j^mr_\tau r_{w}\Psi(F_{w}^{-1}\circ F_\tau^{-1}\circ F_j^{-m}\cdot,F_{w}^{-1}\circ F_\tau^{-1}\circ F_j^{-m}z)\Delta_\mu f(z)d\mu(z).
\end{equation}

     Since $\Psi(\cdot,z)$ is piecewise harmonic, the normal derivatives of $\Psi(\cdot,z)$ are bounded by a constant $c>0$.  So we have 
$$
     |\partial_nr_w\Psi(F_w^{-1}\cdot,F_w^{-1}z)(\cdot)|=|\partial_n\Psi(\cdot,F_w^{-1}z)(F_w^{-1}\cdot)|\leq c.
$$
      Using the above estimate into each summand of $f_3$, we have
         \begin{equation}\label{f2}
     \begin{aligned}
          |d_{i2}f_3(y)|\leq \sum_{0\leq n \leq |\tau|-1}|\int_{F_j^mF_{\tau_1}...F_{\tau_n}K}c\cdot d\mu(z)|\cdot\Vert\Delta_\mu f(z)\Vert_\infty=O(\mu_j^m).
     \end{aligned}
     \end{equation} 
    Using Lemma 3.3, we have an estimate for $f_4$, 
     \begin{equation}\label{f3}
     |d_{i2}f_4(y)|=|\int_{F_j^mF_\tau K}H_i(F_\tau^{-1}F_{j}^{-m}z)\Delta_\mu f(z)d\mu(z)|=O(\mu_j^m\mu_\tau).
     \end{equation}
     
 Finally, Combining \eqref{f1}, \eqref{f2}, \eqref{f3} and Lemma 3.2, we have proved that
$$
      d_{i2}f(y)=d_{i2}f_1(y)+d_{i2}f_3(y)+d_{i2}f_4(y)+d_{i2}h(y)=
      \begin{cases}
      O(\mu_{j}^m), &\text{if $r_j\mu_{j}>|\lambda_{j3}|,$}\\
      O(m\mu_{j}^m), &\text{if $r_j\mu_{j}=|\lambda_{j3}|,$}\\
      O((\lambda_{j3}r_j^{-1})^m),  &\text{if $r_j\mu_{j}<|\lambda_{j3}|.$}
      \end{cases}
$$
Thus we have proved Theorem 1.4. $\Box$

\textit{Remark 1.} The condition $\partial_n f(x)=0$ is necessary. Otherwise, the continuity result in Theorem 1.4 is not true. For example, consider the harmonic function $h=H_2+H_3$, which is a multiple of $h_{12}$, on the Sierpinski gasket, $\mathcal{SG}$, equipped with standard Dirichlet form. It is easy to calculate that $d_{12}h(v_1)=-2$ and $d_{32}h_{12}(F_1^mF_2v_3)=0$ for all $m\geq 0$. Thus $d_{32}h_{12}(F_1^mF_2v_3)$ does not converges to $d_{12}h_{12}(v_1)$, although $F_1^mF_2v_3$ converges to $v_1$, as $m\rightarrow\infty.$ See Figure 3.1 for the values of $h$. 

\begin{figure}[h]
\begin{center}
\includegraphics[width=7cm,totalheight=6cm]{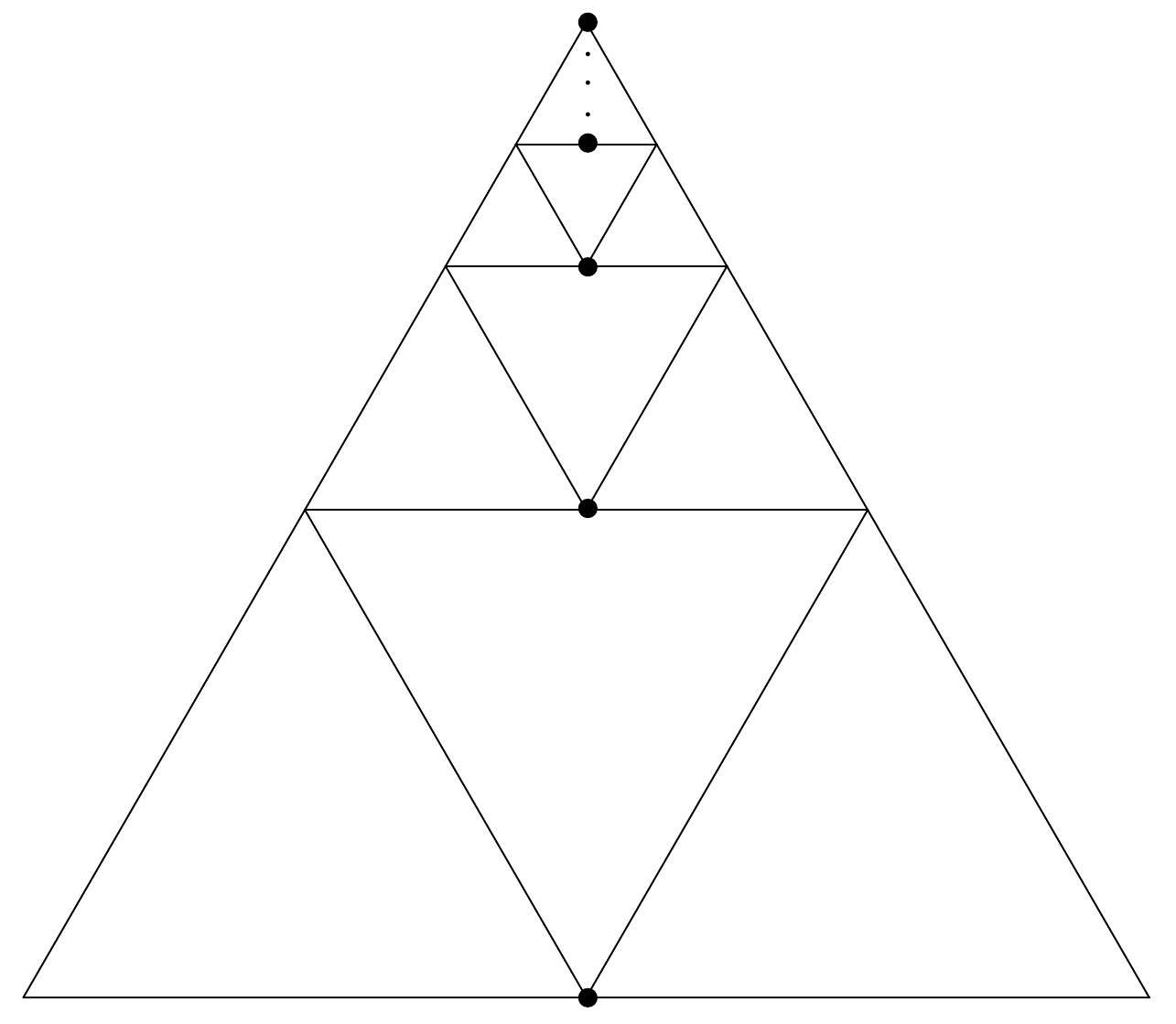}
\setlength{\unitlength}{1cm}
\begin{picture}(0,0) \thicklines
\put(-7.3,-0){$1$}\put(-.3,-0){$1$}\put(-4,-0.2){$4/5$}\put(-4.2,2.7){$12/25$}
\put(-6.0,2.9){$3/5$}\put(-2.0,2.9){$3/5$}\put(-2.8,4.4){$9/25$}\put(-5.3,4.4){$9/25$}\put(-3.82,6){$0$}
\end{picture}
\begin{center}\vspace{0.4cm}
\textbf{Figure 3.1.}
\end{center}
\end{center}
\end{figure}

\textit{Remark 2.} The estimates in Theorem 1.4 are sharp.

In $r_j\mu_j<|\lambda_{j3}|$ case, consider the harmonic function $h_{j3}$. We can find a vertex $z$ in $K$ with nonzero normal derivative, and then  $\partial_n h_{j3}(F_j^m z)=(\lambda_{j3}r_j^{-1})^m\partial_n h_{j3}(z)$.

In  $r_j\mu_j>|\lambda_{j3}|$ case, we take a function $f\in dom(\Delta_\mu)$, satisfying $\Delta_\mu f(x)\equiv1$ on $K$ and $\partial_nf(v_j)=0$. Look at the points $F_j^mv_i, i\neq j$. We have  $\sum_{i\neq j}\partial_nf(F_j^mv_i)=\mu_j^m$, by using the Gauss-Green formula.

As for $r_j\mu_j=|\lambda_{j3}|$ case, we divide the situation into two possible cases, depending on whether we have $\int_{K}a_{j3}(z)d\mu(z)=0$.

In the case that $\int_{K}a_{j3}(z)d\mu(z)\neq 0$  and $\lambda_{j3}>0$,  we still look at the function $f$ with $0$ normal derivate at $v_j$ and $\Delta_\mu f\equiv1$. Fix a vertex $z$ satisfying $\partial_nh_{j3}(z)\neq 0$.  Then
$$   \begin{aligned}
   d_{i2}f(F_j^mz)&=d_{i2}f_1(F_j^mz)+d_{i2}f_2(F_j^mz)+d_{j2}h(F_{j}^mz)\\
                  &=-\sum_{n=0}^{m-1}r_{j}^n\lambda_{j3}^{-n}\mu_{j}^n(\int_K a_{j3}(z)d\mu(z)) d_{i2}h_{j3}(F_j^mz)+O(\mu_j^m)\\
                  &=-mr_j^{-m}\lambda_{j3}^m(\int_K a_{j3}(z)d\mu(z)) d_{i2}h_{j3}(z)+O(\mu_j^m)\\
                  &=O(m\mu_j^m),
   \end{aligned}
$$
where $f_1$, $f_2$ and $h$ are as same as those in the proof of Theorem 1.4.

In the remaining case, (It may happen that $\int_{K}a_{j3}(z)d\mu(z)=0$, for example, if we choose the Sierpinski gasket, $\mathcal{SG}$, equipped with the standard Dirichlet form.), we give another example which looks somewhat complicated.      

\textit{Example 3.6.} Let $\{c_n\}_{n\geq 0}$ be a sequence of numbers with $c_n>0$, $c_n\rightarrow 0$, and $\sum_{n\geq 0}c_n=+\infty.$
Let $\phi$ be a nonnegative continuous function on $K$ satisfying 
\begin{equation}
c_n |\lambda_{j3}|^nr_j^{-n}\leq \min_{F_j^nK\setminus F_{j}^{n+1}K}\phi\leq \max_{F_{j}^nK\setminus F_{j}^{n+1}K}\phi\leq 2c_n|\lambda_{j3}|^n r_j^{-n}.\label{phi}
\end{equation}
Let $g$ be a function on $K$, defined as
\begin{equation}
g(x)=\phi(x)\sum_{n=0}^{\infty}\lambda_{j3}^{-n}r_j^n a_{j3}\circ F_j^{-n}(x).\label{gx}
\end{equation} 
Obviously, $g$ is continuous on $K$, and $g(x)\rightarrow 0$ as $x\rightarrow v_j$, since $\max_{F_j^nK}\phi=o((\lambda_{j3}r_j^{-1})^n)$.
Define
$$
  f(x)=-\int_K(G(x,z)+h_{j2}(x)H_j(z))g(z)d\mu(z).
$$
It is easy to check that $\Delta_\mu f=g$ and $\partial_n f(v_j)=0$.

Let $f=f_1+f_2$ be as those in the proof of Theorem 1.4. Let $y\in F_j^mK$. Following the proof of Theorem 1.4, $\partial_n f_2(y)=O(\mu_j^m)$.

Recalling  the exact formula \eqref{31} for $f_1$, we write
 $$f_1(\cdot)=-h_{j3}(\cdot)\int_{K}\sum_{n=0}^{m-1}r_j^n\lambda_{j3}^{-n}a_{j3}\circ F_j^{-n}g(z)d\mu(z)+R(\cdot),$$ where $R$ denotes the summation of  the remaining terms in \eqref{31}.  
 Following  the arguments in the proof of Theorem 1.4, we have also $\partial_n R(y)=O(\mu_j^m)$.

So we only need to estimate $$\partial_nh_{j3}(y)\int_{K}\sum_{n=0}^{m-1}r_j^n\lambda_{j3}^{-n}a_{j3}\circ F_j^{-n}(z)g(z)d\mu(z).$$ By the proof of Lemma 3.2, $\partial_nh_{j3}(y)=O(\mu_j^m)$.

As for the coefficient $I:=\int_{K}\sum_{n=0}^{m-1}r_j^n\lambda_{j3}^{-n}a_{j3}\circ F_j^{-n}(z)g(z)d\mu(z)$, we write it as
$I=I_1+I_2$, where
$$I_1=\int_{F_j^mK}\sum_{n=0}^{m-1}r_j^n\lambda_{j3}^{-n}a_{j3}\circ F_j^{-n}(z)g(z)d\mu(z)$$ and 
$$I_2=\sum_{l=0}^{m-1}\int_{F_j^lK\setminus F_{j}^{l+1}K}\sum_{n=0}^{l}r_j^n\lambda_{j3}^{-n}a_{j3}\circ F_j^{-n}(z)g(z)d\mu(z).$$
It is easy to verify that $|I_1|=o(r_j^m\lambda_{j3}^{-m}\mu_j^m)=o(1)$, since $g(z)\rightarrow 0$ as $z\rightarrow v_j$.

Taking the expression \eqref{gx} of $g$ into $I_2$, we have
$$I_2=\sum_{l=0}^{m-1}\int_{F_j^lK\setminus F_{j}^{l+1}K}(\sum_{n=0}^{l}r_j^n\lambda_{j3}^{-n}a_{j3}\circ F_j^{-n})^2\phi(z) d\mu(z).$$
Using the estimate \eqref{phi} of $\phi$ on each $F_j^lK\setminus F_j^{l+1}K$, we get
$$I_2\geq \sum_{l=0}^{m-1} c (r_j^l|\lambda_{j3}|^{-l})^2 c_l|\lambda_{j3}|^l r_{j}^{-l}\mu_j^l=c\sum_{n=0}^{m-1}c_n$$ for some constant $c>0$.

Combining all the above estimates, we finally obtain that
$$|\partial_n f(y)|=O(\sum_{n=0}^{m-1}c_n\mu_j^m).$$

Looking at the choice of $\{c_n\}$, we have an estimate of $|\partial_n f(y)|$ which could be very close to $O(m\mu_j^m)$, although it still equals to $o(m\mu_j^m)$.

\section{Boundness and weak continuity of other derivatives}

In this section, we prove Theorem 1.5 and Theorem 1.6. Also, we provide some examples under the proofs.

\textit{Proof of Theorem 1.5.}  (a)  From Proposition 2.1(e), we have $\sum_{l=1}^{N_0}(\beta_{jk})_l=0$, $k\geq 2$. Combining it  with the fact that  $h$ satisfies the  H\"{o}lder estimate that $|h(x)-h(y)|\leq c r_w$ for any $x,y\in F_wK$, with some constant $c>0$, we have
$$
   \begin{aligned}
   &|d_{jk}h(x)|=|r_w^{-1}\beta_{jk}f|_{F_wV_0}|\leq c \text{ for nonjunction vertices, and }\\
   &|d_{j'k}h(x)|=|r_w^{-1}r_j^{-1}\beta_{j'k}f|_{F_wF_jV_0}|\leq c
   \text{ for junction vertices.}
   \end{aligned} 
$$

(b)  The differentiability of $f$ at all vertices is provided by Theorem 4.1 in {[}S5{]}. We now prove the boundness property of $f$. Let $x=F_wv_j$ be a nonjunction vertex. 
We use the notations  in Section 3 that $a_{jk}(z)=d_{jk}\Psi(v_j,z)$. We have
   \begin{equation}
   -\int_K a_{jk}(z)\Delta_\mu f(z)d\mu(z)=\lambda_{jk}^{-1}\beta_{jk}f|_{F_jV_0}-\beta_{jk}f|_{V_0}.\label{ajk}
   \end{equation}
   In fact, it is an immediate result of Lemma 3.4, where we choose $m=1$, and replace $u$ by $f$.
   Scaling  \eqref{ajk} down to $F_wF_j^nK$, $n\geq 0$, we get
\begin{equation}\label{ccc}
\begin{aligned}
&-\int_{F_wF_j^nK}r_j^n\lambda_{jk}^{-n} a_{jk}\circ F_j^{-n}\circ F_w^{-1}(z)\Delta_\mu f(z)d\mu(z)\\=&r_w^{-1}\lambda_{jk}^{-(n+1)}\beta_{jk}f|_{F_wF_j^{n+1}V_0}-r_w^{-1}\lambda_{jk}^{-n}\beta_{jk}f|_{F_wF_j^nV_0}.   
\end{aligned}
\end{equation}
Summing \eqref{ccc} from $n=0$ to $m-1$, we have
\begin{equation}\label{ddd}
   \begin{aligned}
      -\sum_{n=0}^{m-1} \int_{F_w F_j^n K} r_j^n\lambda_{jk}^{-n}a_{jk}\circ F_j^{-n}&\circ F_w^{-1}(z)\Delta_\mu f(z)d\mu(z)=\\&r_w^{-1}\lambda_{jk}^{-m}\beta_{jk}f|_{F_wF_j^mV_0}-r_w^{-1}\beta_{jk}f|_{F_wV_0}.
   \end{aligned}
   \end{equation}
   
   Since $f$ is differentiable at $x$, the limit of the left side of \eqref{ddd} exists as $m\rightarrow\infty$. Moreover, by using the assumption that $r_j\mu_j<|\lambda_{jN_0}|$, it can be bounded as
      \begin{equation}\label{44}
   \begin{aligned}
   |\sum_{n=0}^{\infty} \int_{F_w\circ F_j^n K} r_j^n\lambda_{jk}^{-n}a_{jk}\circ F_j^{-n}\circ F_w^{-1}(z)\Delta_\mu f(z)d\mu(z)|\\\leq\mu_w \sum_{n=0}^{\infty}\lambda_{jk}^{-n}r_j^n\mu_j^n\Vert a_{jk}\Vert_{\infty}\Vert \Delta_\mu f\Vert_{\infty}\leq\mu_w c_1
   \end{aligned}
   \end{equation}
   with some constant $c_1>0$ for all $k\geq 2$.
   
On the other hand, similar to those in the proof of (a) part, by using the H\"{o}lder estimate property of $f$, we also have 
$|r_w^{-1}\beta_{jk}f|_{F_wV_0}|\leq c_2$ for some constant $c_2>0$.
    
Thus
   \begin{equation}\label{45}
      d_{jk}f(x)=-\sum_{n=0}^{\infty} \int_{F_w F_j^n K} r_j^n\lambda_{jk}^{-n}a_{jk}\circ F_j^{-n}\circ F_w^{-1}(z)\Delta_\mu f(z)d\mu(z)+r_w^{-1}\beta_{jk}f|_{F_wV_0}
   \end{equation}
   is uniformly bounded. For the junction vertices, it is also true by using a similar argument. 
   Thus, all derivatives of $f$ are uniformly bounded over all vertices on $K. \Box$
   
 \textit{Proof of Theorem 1.6.}  The proof is analogous to that of Lemma 3.2 and Theorem 1.4, with suitable modifications.  We still give the proof for $x=v_j$, since for other vertices, we could use scaling. 
 
 (a)  For any harmonic function $h$, and any point $y\in F_j^mK\setminus \{v_j\}$, we have the following equality using scaling,
$$
   d_{ik}(h\circ F_j^m)(F_j^{-m}y)=r_j^md_{ik}h(y).
$$
   Since $d_{ik}h_{jl}$ is uniformly bounded by a constant $c>0$ for all $l\geq 3$, as guaranteed by Theorem 1.5, we have
      \begin{equation}\label{eee}
   \begin{aligned}
   |d_{ik}h_{jl}(y)|&=|r_j^{-m}d_{ik}(h_{jl}\circ F_j^m)(F_j^{-m}y)|\\
   &=|r_j^{-m}\lambda_{jl}^md_{ik}h_{jl}(F_j^{-m}y)|\\
   &\leq c r_j^{-m}|\lambda_{jl}|^m\leq c(|\lambda_{j3}|r_j^{-1})^m
   \end{aligned}
   \end{equation}
   for all $y\in F_j^mK\setminus\{v_j\}$, where we use  Proposition 2.2(a). 
   
   On the other hand, Since $h$ assumes $0$ normal derivative at $v_j$, by using Proposition 2.2(c), we could write
$$h=h(v_j)+\sum_{l=3}^{N_0}d_{jl}h(v_j)h_{jl}.$$ Combining this with \eqref{eee}, we have $d_{ik}h(y)=O((\lambda_{j3}r_j^{-1})^m)$ for all $y\in F_j^mK\setminus\{v_j\}$.

(b)  Fix a vertex $y=F_j^mF_\tau v_i$. Formule \eqref{44} and \eqref{45} say that 
\begin{equation}\label{47}
|d_{ik}f(y)-r_j^{-m}r_{\tau}^{-1}\beta_{ik}f|_{F_{j}^{m}F_\tau V_0}|\leq c_1\mu_\tau\mu_j^m=O(\mu_j^m).
\end{equation}  As showed in the proof of Theorem 1.4, we could write $$f=f_1+f_3+f_4+h\text{ on } F_j^m K.$$ 
Thus by using \eqref{47},   
\begin{equation}\label{48}
   \begin{aligned}
   d_{ik}f(y)&=r_j^{-m}r_\tau^{-1}\beta_{ik}f|_{F_j^mF_\tau V_0}+(d_{ik}f(y)-r_j^{-m}r_\tau^{-1}\beta_{ik}f|_{F_j^mF_\tau V_0})\\
   &=r_j^{-m}r_\tau^{-1}\beta_{ik}f|_{F_j^mF_\tau V_0}+O(\mu_j^m)\\
   &=d_{ik}f_1(y)+d_{ik}f_3(y)+r_j^{-m}r_\tau^{-1}\beta_{ik}f_4|_{F_j^mF_\tau V_0}+O((\lambda_{j3}r_j^{-1})^m),
   \end{aligned}
   \end{equation}
   where the last equality follows from the facts that $f_1$ is harmonic on the cell $F_j^mK$, $f_3$ is harmonic on the cell $F_j^mF_{\tau}K$, and the using of (a) part for $d_{ik}h(y)$.
  
  Now we estimate $d_{ik}f_1(y)$, $d_{ik}f_3(y)$ and $r_j^{-m}r_\tau^{-1}\beta_{ik}f_4|_{F_j^mF_\tau V_0}$ separately.
   
   A similar argument as that in the proof of Theorem 1.4 for $f_1$ yields that
   \begin{equation} \label{49}
   d_{ik}f_1(y)=O((\lambda_{j3}r_j^{-1})^m)
   \end{equation}
    for all $y\in F_j^mK\setminus \{v_j\}$.
   
To estimate  $d_{ik}f_3(y)$, we refer to the following scaling.
   \begin{equation}\label{410}
   \begin{aligned}
   &|d_{ik}\Psi(F_{\tau_n}^{-1}...F_{\tau_1}^{-1}\circ F_j^{-m}\cdot,F_{\tau_n}^{-1}...F_{\tau_1}^{-1}\circ F_j^{-m}z)(y)|\\
   =&r_j^{-m}r_{\tau_1}^{-1}...r_{\tau_n}^{-1}|d_{ik}\Psi(\cdot,F_{\tau_n}^{-1}...F_{\tau_1}^{-1}\circ F_j^{-m}z)(F_{\tau_n}^{-1}...F_{\tau_1}^{-1}\circ F_j^{-m}y)|\\
   \leq& c r_j^{-m}r_{\tau_1}^{-1}...r_{\tau_n}^{-1},
   \end{aligned}
   \end{equation}
   for some constant $c>0$ by using Theorem 1.5. 

So by using the  exact expression \eqref{310} of $f_3$,
we have that
\begin{equation}\label{411}
   \begin{aligned}
  d_{ik}f_3(y)
  =&-\sum_{0\leq n\leq |\tau|-1}\int_{K} r_j^mr_{\tau_1}...r_{\tau_n}\\&d_{ik}\Psi(F_{\tau_n}^{-1}...F_{\tau_1}^{-1}\circ F_j^{-m}\cdot,F_{\tau_n}^{-1}...F_{\tau_1}^{-1}\circ F_j^{-m}z)(y)\Delta_\mu f(z)d\mu(z)\\
  &=O(\mu_j^m).
   \end{aligned}
   \end{equation}

As for $r_j^{-m}r_\tau^{-1}\beta_{ik} f_4|_{F_j^mF_\tau V_0}$, we observe that  it equals $0$ since $f_4$ takes zero values on the boundary of $F_j^mF_\tau K$. 

Combining the above observation  with \eqref{48}, \eqref{49}, \eqref{410} and \eqref{411}, and the fact that $\mu_j<r_{j}^{-1}|\lambda_{jN_0}|<r_j^{-1}|\lambda_{j3}|$, we have finally proved the (b) part of the theorem. $\Box$

\textit{Remark 1.} Suppose $\sharp V_0=3$ and all structures have full $D3$ symmetry. Theorem 1.5(b) and Theorem 1.6(b) are still valid without the hypothesis $r_j\mu_{j}<|\lambda_{j3}|$ (in this case, $N_0=3$), if we additional assume that $g=\Delta_\mu f$ satisfies the H\"{o}lder condition that
$$|g(x)-g(y)|\leq c\gamma ^m$$ for all $x,y$ belonging to the same $m$-cells, for some constant $\gamma$ satisfying 
\begin{equation}\label{ttt}
r_j\mu_j\gamma<|\lambda_{j3}|,
\end{equation}
 for all $j$. 

The key observation is that $a_{j3}$ is skew-symmetry with respect to the point $v_j$, which yields that in \eqref{44}, each term in the summation could be rewrote as, 
$$\int_{F_wF_j^nK}r_j^n\lambda_{j3}^{-n}a_{j3}\circ F_j^{-n}\circ F_w^{-1}(z)(\Delta_\mu f(z)-\Delta_\mu f(x))d\mu(z),$$
and this is estimated by a multiple of $\mu_{w}r_j^n|\lambda_{j3}|^{-n}\mu_j^n\gamma^n$. 
Since $r_j\mu_j\gamma<|\lambda_{j3}|$, this guarantees the convergence of \eqref{44}. The existence of the derivatives also holds, which was proved in {[}S5{]}, due to the same reason.

\textit{Example 4.1.}  (1) The Sierpinski gasket, which has all $r_j=3/5$, $\mu_j=1/3$, $\lambda_{j3}=1/5$ in the $D3$ symmetry case. Hence $r_j\mu_j=\lambda_{j3}$ for all $j$.

(2) The hexagasket, which can be generated by 6 mappings with simultaneously rotate and contract by a ratio of $1/3$ in the plane. In this case, we take all $r_j=3/7$, $\mu_j=1/6$ and $\lambda_{j3}=1/7$, thus the condition $r_j\mu_j<|\lambda_{j3}|$ holds. See Figure 4.1 for the first two level graphs that approximate the hexagasket. 

(3) The level $3$ Sierpinski gasket, $\mathcal{SG}_3$, obtained by taking 6 contractive mappings of ratios $1/3$, as shown in Figure 4.2. All $r_j=7/15$, $\mu_j=1/6$ and $\lambda_{j3}=1/15$. Thus the condition $r_j\mu_j<|\lambda_{j3}|$ does not hold.

 Please find the detail information of these examples in the book {[}S8{]}. If $\Delta_\mu f\in dom(\Delta_\mu)$ then \eqref{ttt} holds with $\gamma=r_j$ as showed in \eqref{holder}. This holds in examples (1) and (3) above. Thus the conclusions in Theorem 1.5 and 1.6 are valid for these fractals.
 
 \begin{figure}[h]
\begin{center}
\includegraphics[width=10cm,totalheight=5cm]{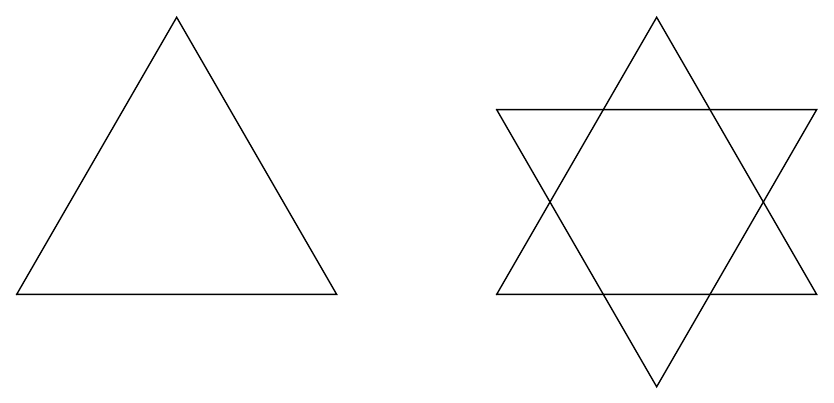}
\setlength{\unitlength}{1cm}
\begin{center}\vspace{0.4cm}
\textbf{Figure 4.1.} The first 2 graphs that approximate the hexagasket.
\end{center}
\end{center}
\end{figure}

 \begin{figure}[h]
\begin{center}
\includegraphics[width=5.6cm,totalheight=5cm]{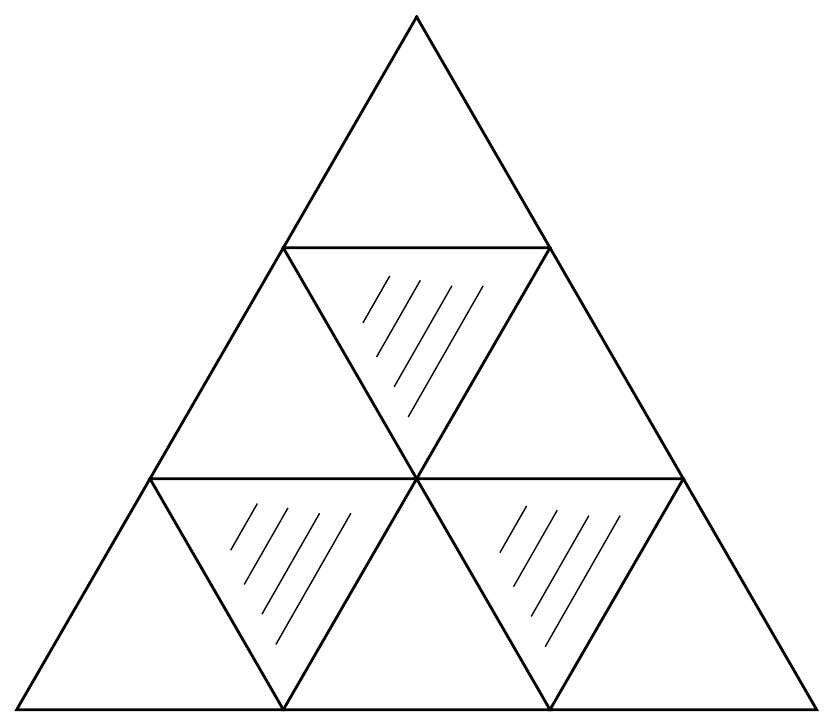}
\setlength{\unitlength}{1cm}
\begin{center}\vspace{0.4cm}
\textbf{Figure 4.2.} The first  graph that approximate $\textit{SG}_3$.
\end{center}
\end{center}
\end{figure}

\textit{Remark 2.}  The condition $d_{j2} f(x)=0$ in Theorem 1.6 could not be replaced by $d_{jk}f(x)=0$, although it looks more "reasonable". For example, look at the Sierpinski gasket, $\mathcal{SG}$, equipped with the standard Dirichlet form. We consider the harmonic function $h=H_2+H_3$, which is a multiple of $h_{12}$. It is easy to calculate that $d_{12}h(v_1)=-2$, $d_{13}h(v_1)=0$, and $d_{13}h(F_1^mv_2)=1/3$ for all $m\geq 1$. Thus $d_{13}h(F_1^{m}v_2)$ does not converges to $d_{13}h(v_1)$, although $F_1^mv_2$ converges to $v_1$, as $m\rightarrow\infty$. See Figure 4.3 for the values of $h$.

\begin{figure}[h]
\begin{center}
\includegraphics[width=7cm,totalheight=6cm]{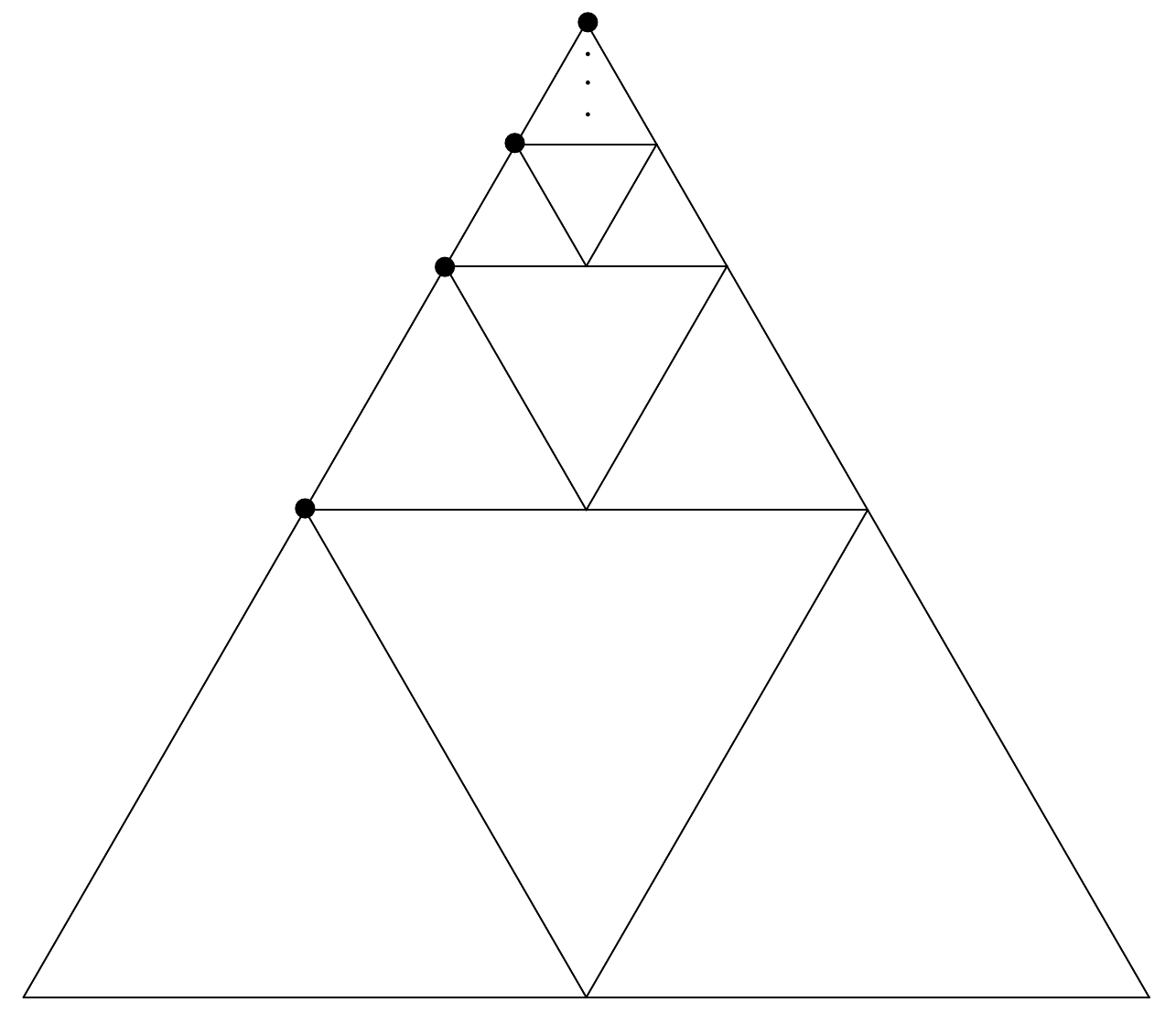}
\setlength{\unitlength}{1cm}
\begin{picture}(0,0) \thicklines
\put(-7.3,-0){$1$}\put(-.3,-0){$1$}\put(-4,-0.2){$4/5$}\put(-4.2,2.7){$12/25$}
\put(-6.0,2.9){$3/5$}\put(-2.0,2.9){$3/5$}\put(-2.8,4.4){$9/25$}\put(-5.3,4.4){$9/25$}\put(-3.82,6){$0$}
\end{picture}
\begin{center}\vspace{0.4cm}
\textbf{Figure 4.3.}
\end{center}
\end{center}
\end{figure}

\textit{Remark 3.}  As we know, the assumption $r_i\mu_i<|\lambda_{iN_0}|$ in Theorem 1.5(b) is only a sufficient condition which guarantees the existence of all derivatives of $f$. It could be relaxed as stated in Remark 1 in the $D3$ symmetry case. One may ask a question that: Whether does Theorem 1.5(b) still hold as long as $f\in dom(\Delta_\mu)$ and $f$ is differentiable at all vertex? We will give an example  to illustrate that this is not true.

\textit{ Example 4.2.} Consider the Sierpinski gasket, $\mathcal{SG}$, equipped with the standard Dirichlet form and the standard self-similar measure. So all $r_i=3/5, \mu_i=1/3$. 

First, we define a sequence of functions $g_l, l\geq 0$, satisfying
$$-\Delta_\mu g_l(x)=\sum_{n=0}^{l} a_{33}\circ F_3^{-n}(x)$$ with the Dirichlet boundary condition, i.e.,
$g_l|_{V_0}=0.$ Here each term in the summation has the understanding that $a_{33}F_3^{-n}(x)$ is zero unless $x$ belongs to $F_3^n\mathcal{SG}$.

It is easy to observe that $\Vert\Delta_\mu g_{l}\Vert_{\infty}$ is uniformly bounded and $$d_{33}g_l(v_3)>(l+1)c>0,$$ for all $l$ with some constant $c>0$. In fact, by using \eqref{45},
\begin{equation}\label{412}
\begin{aligned}
d_{33}g_l(v_3)=&-\sum_{m=0}^{\infty}\int_{F_3^m\mathcal{SG}}\lambda_{33}^{-m}r_3^{m}a_{33}\circ F_3^{-m}(z)\Delta_\mu g_l(z)d\mu(z)\\
=&\sum_{n=0}^l\sum_{m=0}^{\infty}\int_\mathcal{SG}\lambda_{33}^{-m}r_3^m a_{33}\circ F_3^{-m}(z)a_{33}\circ F_3^{-n}(z)d\mu(z)\\
\geq &\sum_{n=0}^l\sum_{m=n}^{\infty}\int_{F_3^{m}\mathcal{SG}}\lambda_{33}^{-m}r_3^m a_{33}\circ F_3^{-m}(z)a_{33}\circ F_3^{-n}(z)d\mu(z)\\
=&\sum_{n=0}^l\sum_{m=n}^{\infty}\int_{F_3^{m-n}\mathcal{SG}}\lambda_{33}^{-m}r_3^m\mu_3^n a_{33}\circ F_3^{-m+n}(z)a_{33}(z)d\mu(z)\\
=&\sum_{n=0}^l\sum_{m=0}^\infty\int_{F_3^m\mathcal{SG}}\lambda_{33}^{-m} r_3^m a_{33}\circ F_3^{-m}(z)a_{33}(z)d\mu(z)=(l+1)d_{33}g_0(v_3)>0,
\end{aligned}
\end{equation}
noticing that $a_{33}$ is skew-summery with respect to $v_3$.

Now we define a function $g$, which is the solution of the following Dirichlet problem,
$$
    \begin{cases}
    \Delta_\mu g(x)=\sum_{l=0}^{\infty} 3^{-l}\Delta_\mu g_{3^{3l}}\circ F_1^{-1}\circ F_2^{-l}(x),\\
        g|_{V_0}=0.
    \end{cases}
$$
See Figure 4.4 to find the support of $ \Delta_\mu g(x)$. 

\begin{figure}[h]
\begin{center}
\includegraphics[width=7cm,totalheight=6cm]{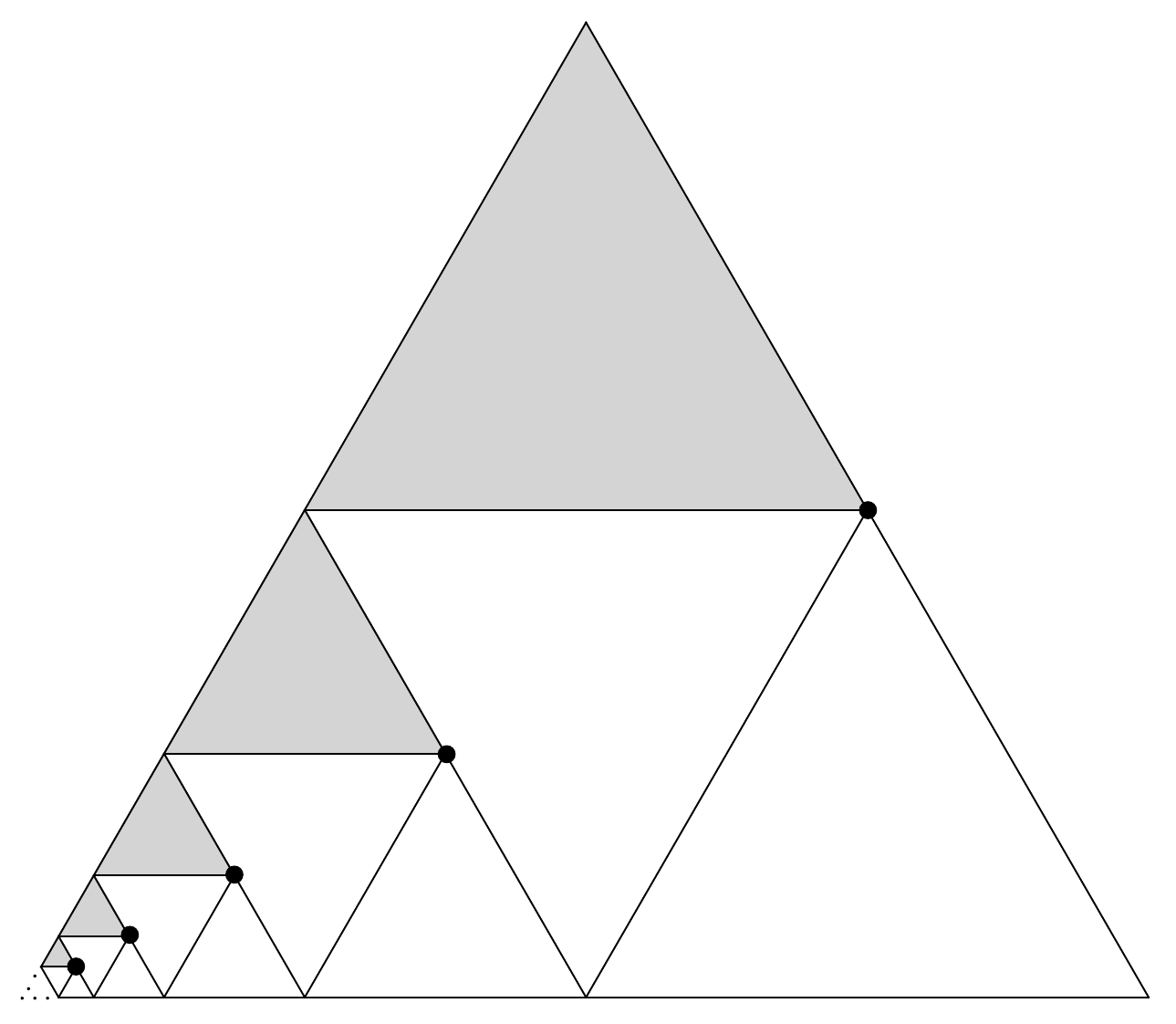}
\setlength{\unitlength}{1cm}
\begin{picture}(0,0) \thicklines
\put(-7.4,-0){$v_2$}\put(-0.3,-0){$v_3$}\put(-3.9,5.95){$v_1$}
\end{picture}
\begin{center}\vspace{0.4cm}
\textbf{Figure 4.4.}
\end{center}
\end{center}
\end{figure}

Next we estimate the tangential derivatives of $g$ at the vertices $F_2^lF_1v_3$. By using \eqref{45}, \eqref{412}, we have
$$
    \begin{aligned}
    d_{33}g(F_2^lF_1v_3)=&-\sum_{m=0}^{\infty}\int_{F_2^lF_1F_3^m\mathcal{SG}} r_3^m\lambda_{33}^{-m}a_{33}\circ F_3^{-m}\circ F_1^{-1}\circ F_2^{-l}(z) \Delta_\mu g(z)d\mu(z) \\
    &+r_2^{-l}r_1^{-1}\beta_{33}g|_{F_2^lF_1V_0}\\
         =&-\sum_{m=0}^\infty\int_{F_3^m\mathcal{SG}}3^{-l}\mu_2^l\mu_1r_3^m\lambda_{33}^{-m}a_{33}\circ F_3^{-m}(z)\Delta_\mu g_{3^{3l}}(z)d\mu(z) \\
    &+r_2^{-l}r_1^{-1}\beta_{33}g|_{F_2^lF_1V_0}\\
    =&-3^{-2l-1}\sum_{m=0}^{\infty}\int_{F_3^m\mathcal{SG}} r_3^m\lambda_{33}^{-m}a_{33}\circ F_3^{-m}(z) \Delta_\mu g_{3^{3l}}(z)d\mu(z)+O(1)\\
    =&3^{-2l-1}d_{33}g_{3^{3l}}(v_3)+O(1)\geq c3^{-2l-1}(3^{3l}+1)+O(1).
    \end{aligned}
$$

Thus we have proved that $\{d_{33}g(F_2^l F_1v_3)\}_{l\geq 0}$ is unbounded, although we have $g\in dom(\Delta_\mu)$ and is differentiable at all vertices. (In fact, $\Delta_\mu g$ satisfies the H\"{o}lder continuous condition in a neighborhood at any vertex.)

\section{On the weak tangent}

In this section, we focus on the concept weak tangent.

Let $f$ be a function which is differentiable at a vertex $x$. Then the weak tangent of order one of $f$ at $x$, denoted as $T_1^x(f)$, is the harmonic function on $U_0(x)$ with the same value and the same gradient as $f$ at $x$. Let $h_m$ be the harmonic function assumes the same values as $f$ at the boundary of $U_{m}(x)$, extended to be harmonic on $U_0(x)$. Theorem 3.11 in {[}S5{]} says that $h_m$ converges to $T_1^x(f)$ uniformly on $U_0(x)$ as $m$ goes to infinity.  However, the following example will show that this is not true. 

\textit{Example 5.1.} Consider the Sierpinski gasket $\mathcal{SG}$, equipped with a self-similar Dirichlet form which only has a single bilateral symmetry, as described in Example 2.4.

Define a function $f$ on $\mathcal{SG}$ as following. We assume
$$
   \begin{cases}
   f(F_2F_3^mv_j)=\eta^m(\alpha_{32})_j\quad\text{for $j=1,2$ and $m\geq 0$,} \\
   f(v_1)=0, f(v_3)=0, f(F_1v_3)=0,
   \end{cases}
$$
where $\eta$ is a constant such that $|\lambda_{23}|=|\lambda_{33}|<\eta<\lambda_{22}=\lambda_{32}$. As for the values of $f$ at other points, we take harmonic extension.  

Choose $x=F_2v_3$, it is easy to check that 
$$d_{22}f(x)=d_{23}f(x)=d_{32}f(x)=d_{33}f(x)=0.$$ Thus $f$ is differentiable at $x$ and $T_1^x(f)\equiv 0$ on $U_0(x)$.

On the other hand, using the bilateral symmetry, we could obtain that
$$
   h_m(x)=\frac{\sum_{y\sim_{m+1} x}c_{xy}f(y)}{\sum_{y\sim_{m+1} x}c_{xy}}
         =\eta^m\frac{\sum_{y\sim_1 x}c_{xy}f(y)}{\sum_{y\sim_1 x}c_{xy}}=\eta^mh_0(x),    
$$
which results that 
$$d_{23}h_m(x)=r_3^{-1}\lambda_{23}^{-m}(\beta_{23})_2 h_m(x)=r_3^{-1}\lambda_{23}^{-m}\eta^m(\beta_{23})_2h_0(x).$$
Thus $d_{23}h_m(x)\rightarrow\infty$ as $m\rightarrow\infty$ since $|\lambda_{23}|<\eta$ and $(\beta_{23})_2\neq 0$ as showed in Example 2.4. So we have
$$\beta_{23}h_m|_{F_3V_0}\rightarrow\infty\text{ as } m\rightarrow\infty,$$
which means $\Vert h_m\Vert_\infty\rightarrow\infty\text{ as }m\rightarrow\infty.$
Hence $h_m$ does not converge to $T_1^x(f)$ as $m\rightarrow\infty$.

We need some extra assumption to make the theorem holds. 

\textbf{Theorem 5.2.} \textit{Suppose one of the condition in Proposition 2.3 holds. Then for any $f$ differentiable at $x$, $h_m$ converges to $T_1^x(f)$ uniformly. }

\textit{Proof.} The proof is essential the same as that of Theorem 3.11 in {[}S5{]}, where the condition $(\beta_{jk})_j=0$ may be misapplied. we omit it here. $\Box$

As pointed out under the proof of Proposition 2.3, in the D3 symmetry case, the assumption in Theorem 5.2 holds automatically.

\textbf{Theorem 5.3.} \textit{Suppose
\begin{equation}
r_j\max_{1\leq j\leq N_0}\mu_j<|\lambda_{jN_0}|\label{diff}
\end{equation} 
for every $j$.
Then for any $f\in dom(\Delta_\mu)$, for any vertex $x$,
$h_m$ converges to $T_1^x(f)$ uniformly. }

\textit{Proof.} Condition \eqref{diff} guarantees the differentiability of $f$ at $x$ by using Theorem 4.1 in {[}S5{]}. 

For a nonjunction point $x$, we have
$$d_{jk}f(x)=\lim_{m\to \infty}d_{jk}h_m(x),
$$
since on the right side of \eqref{nonj} we may replace $f$ by $h_m$ and $h_m$ is harmonic on $U_0(x)$. In particular, this also shows the limit exists. We have $h_m(x)=f(x)$ for all $m$ since $x$ is a boundary point of $U_m(x)$. On the other hand, there is an estimate for harmonic functions, $|h(y)|\leq c(|h(x)|+\Vert dh(x)\Vert)$ uniformly for $y\in U_0(x)$, which is a result of Proposition 2.2(c). Using this estimate for $h_m-T_1^x(f)$, we obtain that $h_m$ converges uniformly on $U_0(x)$ to $T_1^x(f)$.

       If $x$ is a junction point, i.e., $x=F_wF_jv_{j'}$ for all $j\in J(x)$, we no longer have $x$ as a boundary point of $U_m(x)$. We have to give an estimate of $h_m(x)-f(x)$. If we only have the compatibility condition, we only have $h_m(x)-f(x)=o(\lambda_{j'2}^m)$. With the assumption \eqref{diff}, we can say more. 
       
       Let $\psi^{m}_x$ denote the tent function, the piecewise harmonic function in $\mathcal{H}_m$ which takes value $1$ at $x$ and $0$ at all other vertices in $V_{m}$.  
       
   By using the pointwise formula for $\Delta_\mu f$ at $x$, we have
    \begin{equation}\label{abc}
               \begin{aligned}
                  \Delta_\mu f(x)&=\lim_{m\to\infty}\frac{\sum_{\sim_{m}}c_{xy}(f(y)-f(x))}{\int_K\psi^{m}_xd\mu}\\
                  &=\lim_{m\to\infty}\frac{\sum_{j\in J(x)}r_w^{-1}r_j^{-1}\lambda_{j'2}^{-m}\beta_{j'2}f|_{F_wF_jF_{j'}^mV_0}}{\int_K\psi^{m+|w|+1}_{x}d\mu}\\
                  &=\lim_{m\to\infty}\frac{\sum_{j\in J(x)}r_w^{-1}r_j^{-1}\lambda_{j'2}^{-m}(\beta_{j'2})_{j'}(f(x)-h_m(x))}{\int_K\psi^{m+|w|+1}_xd\mu},
               \end{aligned}
            \end{equation}
            in which we use the compatibility condition
       $$
                   \sum_{j\in J(x)}r_w^{-1}r_j^{-1}\lambda_{j'2}^{-m}\beta_{j'2}h_m(x)|_{F_wF_jF_{j'}^mV_0}=0,
         $$
            since $h_m$ is harmonic.

       The integral in \eqref{abc} can be calculated, 
       $$
               \int_K\psi^{m+|w|+1}_xd\mu=\sum_{j\in J(x)}\mu_w \mu_j\mu_{j'}^m\int_KH_{j'}d\mu
         $$
            where $H_j$ denote the harmonic function taking $1$ at $v_j$ and $0$ at other points of $V_0$. Thus the integral converges to zero with the rate $(\mu_{J(x)})^m$, where $\mu_{J(x)}=max_{j\in J(x)}\mu_{j'}$.  Denote $r_{J(x)}=\min_{j\in J(x)}\{r_{j'}\}$, we then have
 $$
            f(x)-h_m(x)=O((r_{J(x)}\mu_{J(x)})^m)
   $$    
   from the convergence of \eqref{abc}.
   
   Combining this estimate with the assumption \eqref{diff}, we get
  $$
            f(x)-h_m(x)=o(\lambda_{j'k}),
   $$   for all $j'$.      
            So we have the following equation as the nonjunction case,
    $$
            \begin{aligned}
            d_{j'k}f(x)=&\lim_{m\to\infty}r_w^{-1}r_j^{-1}\lambda_{j'k}^{-m}\beta_{j'k}h_m|_{F_wF_jF_{j'}^mV_0}\\
            &+\lim_{m\to\infty}r_w^{-1}r_j^{-1}\lambda_{j'k}^{-m}(\beta_{j'k})_{j'}(f(x)-h_m(x))\\
                               =&\lim_{m\to\infty}d_{j'k}h_m(x).
            \end{aligned}
   $$
   
   Using a similar argument as the nonjunction case, we also obtain that $h_m$ converges uniformly on $U_0(x)$ to $T_1^x(f)$. $\Box$

At last, we will give an example which could serve as a counter-example of Conjecture 6.7 in {[}S5{]} on weak tangents of higher order. 

\textit{Example 5.4.}  For the Sierpinski gasket $\mathcal{SG}$, we assume all the structures satisfy the $D3$ symmetry.  In this case, all $r_j=3/5$ and $\mu_j=1/3$, $\rho=1/5$. Define a function $f\in dom(\Delta_\mu)$ which satisfies 
\begin{equation}
    \begin{cases}
    \Delta_\mu f=\sum_{m=0}^{\infty} \eta^m\psi_{F_1^mF_2v_3}^{m+1},\\
        f(v_1)=0,df(v_1)=0,\\
    \end{cases}
    \end{equation}
    where $r<\eta<1$, $\psi_x^m$ is a piecewise harmonic spline in $\mathcal{H}_m$ satisfying $\psi_x^m(y)=\delta_{xy}$ for $y\in V_m$.
    One can easily verify that $d\Delta_\mu f(v_1)=0$. We will show that $f$ does not have a weak tangent at $v_1$ of order $2$. 

In fact, by using the Gauss-Green formula, we have
$$
    f(v_2)+f(v_3)=\int_{K}H_1(x)\Delta_\mu f(x)d\mu(x),
$$
where $H_1$ is the harmonic function satisfying $H_1(v_j)=\delta_{1j}$.  Using scaling, we then have
\begin{equation}
    \begin{aligned}
    f(F_1^mv_2)+f(F_1^mv_3)&=\rho^m\int_{K} H_1 (\Delta_\mu f)\circ F_1^{m}(x)d\mu(x)\\
    &=\rho^m\eta^m(f(v_2)+f(v_3)).\label{mf}
    \end{aligned}
\end{equation}

But from the proof of Lemma 6.2 in {[}S5{]}, for any 2-harmonic function $h$, there exist constants $a,b,c\in\mathbb{R}$ such that
    \begin{equation}
    h(F_1^mv_2)+h(F_1^mv_3)=ar^m+b\rho^m+c(r\rho)^m.\label{2h}
    \end{equation}
   Combining  $\eqref{mf}$ and $\eqref{2h}$, we could claim that it is impossible to have any $2$-harmonic function $h$ satisfying $\eqref{wt1}$, where $n$ is replaced with $2$, since $r<\eta<1$. Thus $f$ does not have a weak tangent at $v_1$ of order $2$.
   
   Before the end of this section, we would like to pose a problem that should be considered.  The Hypothesis 1.1 requires the harmonic structure to be nondegenerate, i.e., all the transformation matrices are nonsingular. This excludes some important fractals such as the Vicsek set. Consider a square with corners $\{v_1,v_2,v_3,v_4\}$ and center $v_5$. Let $F_j$ be contractive mappings  with ratio $1/3$ and fixed point $v_j$. The invariant set of this i.f.s. is called the Vicsek set, denoted by $\mathcal{V}$. Then $N=5$, $N_0=4$ with $V_0=\{v_1,v_2,v_3,v_4\}$. See Figure 5.1 for the second level graph of $\mathcal{V}$. This fractal has $D4$ symmetry. Equip $\mathcal{V}$ with the standard Dirichlet form and standard measure. Then all $r_j=1/3$, $\mu_j=1/5$, and all the transformation matrices $M_j$ are permutations of $M_1$ which is
   $$\left(
  \begin{array}{cccc}
1 & 0 & 0 & 0  \\
\frac{3}{4} & \frac{1}{12} &  \frac{1}{12} &  \frac{1}{12}\\
\frac{1}{2} &  \frac{1}{6} & \frac{1}{6}  & \frac{1}{6}  \\
\frac{3}{4} & \frac{1}{12} & \frac{1}{12} &  \frac{1}{12}  \\
  \end{array}
\right).$$ 
It is easy to calculate that $\lambda_{j2}=1/3, \lambda_{j3}=\lambda_{j4}=0$. Thus this harmonic structure of $\mathcal{V}$ is degenerate. Is there a satisfactory theory of derivatives or gradients on $\mathcal{V}$? Or even on other fractals in degenerate case?

  \begin{figure}[h]
\begin{center}
\includegraphics[width=5cm,totalheight=5cm]{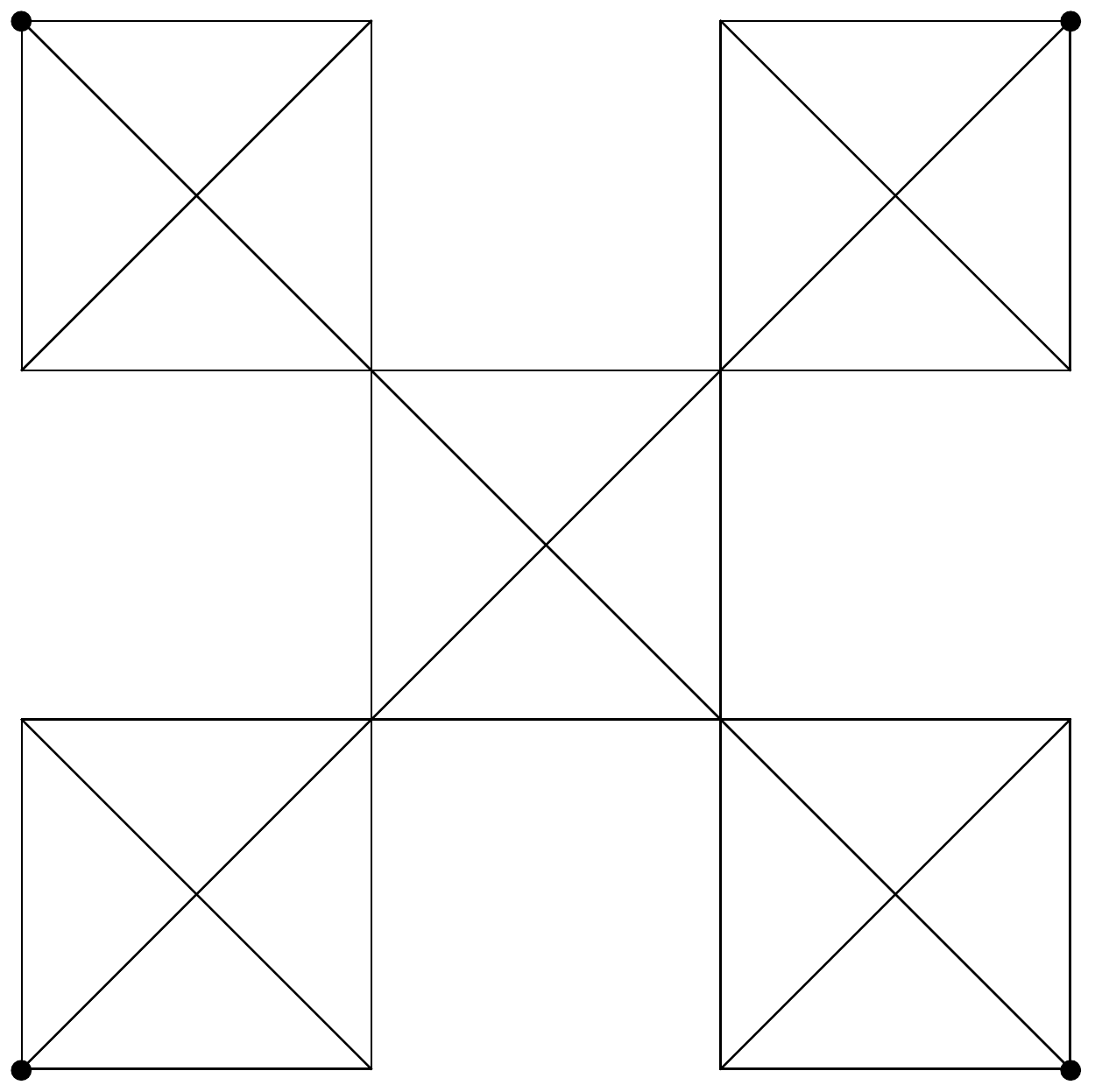}
\setlength{\unitlength}{1cm}
\begin{picture}(0,0) \thicklines
\put(-0.2,-0){$v_3$}\put(-5.6,4.9){$v_1$}\put(-5.6,-0){$v_2$}\put(-0.2,4.9){$v_4$}
\end{picture}
\begin{center}\vspace{0.4cm}
\textbf{Figure 5.1.} The second level graph of $\mathcal{V}$.
\end{center}
\end{center}
\end{figure}

\end{document}